\documentclass[12pt]{article}

\usepackage{amsmath,amssymb}

\def\med{\medbreak\noindent}
\def\sms{\smallskip}
\def\ms{\medskip}

\def\today{\noindent\number\day
\space\ifcase\month\or
  January\or February\or March\or April\or May\or June\or
  July\or August\or September\or October\or November\or December\fi
  \space\number\year}

\def\bP {{\mathbb P}} \def\bE {{\mathbb E}} \def\bQ {{\mathbb Q}}
\def\bR {{\mathbb R}}  \def\bZ {{\mathbb Z}}

  \def\sC {\mathcal{C}}
\def\sD {\mathcal{D}} \def\sE {\mathcal{E}} \def\sF {\mathcal{F}}
\def\sG {\mathcal{G}}  
  \def\sL {\mathcal{L}}

\def\sY {\mathcal{Y}} 

\def\ol{\overline}

 \def\Lam {\Lambda} \def\Gam{\Gamma}

  \def\eps{\varepsilon}
\def\om{\omega }
 
\def\vp{\varphi}

\def\pd {\partial}
\def\q{\quad} \def\qq{\qquad}
\def\dint{\int\kern-.6em\int}
\def\grad{\nabla}

\def\Osc{\mathop{{\rm Osc}}}

\def \fract#1#2{{\textstyle \frac{#1}{#2}}}

\def \half {{\textstyle\frac12}}

\def \qed {\hfill$\square$\par}

\def\=d{{\,\buildrel (d) \over =\,}}
\def\a.s.{{\buildrel a.s. \over \longrightarrow}}

\def\wt{\widetilde}

\newtheorem{theorem}{Theorem}[section]
\newtheorem{corollary}[theorem]{Corollary}
\newtheorem{lemma}[theorem]{Lemma}
\newtheorem{proposition}[theorem]{Proposition}

\newtheorem{assumption}[theorem]{Assumption}

\numberwithin{equation}{section}

\def\vp{{\varphi}}
\def\q{\quad}
\def\qq{{\qquad}}

\def\pd{\partial}

\def\wt{\widetilde}

\def\eps{\varepsilon}

\def\phi{\varphi}

\def\om{\omega}
\def\ol{\overline}
\def\fract{\textstyle \frac}

\def\proof{{\medskip\noindent {\sl Proof. }}}

\def\qed{{\hfill $\square$ \bigskip}}

\def\Osc{\mathop {\rm Osc\, }}

  \def\sC {{\cal C}}
\def\sD {{\cal D}} \def\sE {{\cal E}} \def\sF {{\cal F}}
\def\sG {{\cal G}}  
  \def\sL {{\cal L}}

\def\sY {{\cal Y}}

 \def\bE {{\mathbb E}}

\def\bP {{\mathbb P}} \def\bQ {{\mathbb Q}} \def\bR {{\mathbb R}}

 \def\bZ {{\mathbb Z}}

\def\half{{\textstyle \frac12}}
\def\tfrac{\textstyle \frac}
\def\ignore#1{}

\def\ms{\medskip}
\def\sms{\smallskip}
\def\med{\medskip\noindent}
\def\sm{\smallskip\noindent}
\def\Ci{\sC_\infty}  
\def\Gam{\Gamma}

\begin{document}

\title{\bf Parabolic Harnack Inequality and
 Local Limit Theorem for Percolation Clusters}

\author{
M. T. Barlow\footnote{Research partially supported by NSERC
(Canada), and EPSRC(UK) grant EP/E004245/1} { }and  
B. M. Hambly\footnote{Research supported by EPSRC grant EP/E004245/1}}

\maketitle


\begin{abstract}
We consider the random walk on supercritical percolation clusters in
$\bZ^d$. Previous papers have obtained Gaussian heat kernel bounds,
and a.s. invariance principles for this process. We show how this
information leads to a parabolic Harnack inequality, a local limit
theorem and estimates on the Green's function.


\end{abstract}

\section{Introduction}

We begin by recalling the definition of bond percolation on $\bZ^d$:
for background on percolation see \cite{Grim}. 
We work on the Euclidean lattice $(\bZ^d, \bE_d)$, where $d\ge 2$ and 
$\bE_d =\big\{ \{x,y\}: |x-y|=1\big\}$. Let
$\Omega= \{0,1\}^{\bE_d}$, $p\in [0,1]$, and $\bP=\bP_p$
be the probability measure on $\Omega$ which makes
$\om(e)$, $e \in \bE_d$ 
i.i.d.~Bernoulli r.v., with $\bP(\om(e)=1)=p$.
Edges $e$ with $\om(e)=1$ are called {\sl open}
and the {\sl open cluster} $\sC(x)$ containing $x$ is the set of $y$ such that
$x \leftrightarrow y$, that is $x$ and $y$ are connected by an open path. 
It is well known that there exists $p_c\in (0,1)$ such that
when $p>p_c$ there is a unique infinite open cluster, which we denote
$\Ci=\Ci(\om)$.

Let $X =(X_n, n \in \bZ_+, P^x_\om, x\in \Ci)$ be the
simple random walk (SRW) on $\Ci$.  At each time step,
starting from a point $x$, the process $X$ jumps 
along one of the open edges $e$ containing $x$, with each
edge chosen with equal probability.
If we write $\mu_{xy}(\om)=1$ if $\{x,y\}$ is an open edge
and 0 otherwise, and set $\mu_x =\sum_y \mu_{xy}$, 
then $X$ has transition probabilities
\begin{equation}\label{PXdef}
 P_X(x,y) = \frac{\mu_{xy}}{\mu_x}. 
\end{equation}
We define the transition density of $X$ by
\begin{equation}
  p^\om_n(x,y) = \frac{P^x_\om(X_n=y)}{\mu_y}.
\end{equation}
This random walk on the cluster $\Ci$ was called by 
De Gennes in \cite{DG} `the ant in the labyrinth'. 

Subsequently slightly different walks have been considered: 
the walk above is called the `myopic ant', while there is also a
version called the `blind ant'.
See \cite{Ma}, or Section 5 below for a precise definition.

There has recently been significant progress in the study of this
process, and the closely related continuous time random walk
$Y=(Y_t, t \in [0,\infty), \tilde P^x, x \in \Ci)$, with
generator
$$ \sL f(x)= \sum_y \frac{\mu_{xy}}{\mu_x} (f(y)-f(x)). $$
We write  
\begin{equation}
  q^\om_t(x,y) = \frac{\tilde P^x_\om(Y_t=y)}{\mu_y}
\end{equation}
for the transition densities of $Y$.
Mathieu and Remy in \cite{MR} obtained a.s. upper 
bounds on $\sup_y q^\om_t(x,y)$, and these were extended in \cite{B1}
to full Gaussian-type upper and lower bounds  -- see 
\cite[Theorem 1.1]{B1}.
A quenched or a.s.  invariance principle for $X$ was then obtained in
\cite{SS, BB, MP}: an averaged, or annealed invariance principle
had been proved many years previously in \cite{DFGW}.

The main result in this paper is that as well
as the invariance principle, one also has a local limit
theorem for $p^\om_n(x,y)$ and $q^\om_t(x,y)$. 
(See \cite{Fel}, XV.5 for the
classical local limit theorem for lattice r.v.) For $D>0$ write 
$$ k^{(D)}_t(x) = (2\pi t D)^{-d/2} e^{-|x|^2/2Dt} $$
for the Gaussian heat kernel with diffusion constant $D$.

\begin{theorem}\label{thm:llt-pi} Let $X$ be either the
`myopic' or the `blind' ant random walk on $\Ci$. Let $T>0$. 
Let $g^\om_n: \bR^d \to \Ci(\om)$ be defined so that 
$g^\om_n(x)$ is a 
closest point in $\Ci(\om)$ to $\sqrt{n}x$.
Then there exist constants $a$, $D$ (depending only on $d$ and $p$,
and whether $X$ is the blind or myopic ant walk) such that
$\bP$-a.s. on the event $\{0 \in \Ci\}$,
\begin{equation} \label{llt-pvi}
  \lim_{n\to\infty}  
\sup_{x\in \bR^d} \sup_{t \ge T} 
 \Big|n^{d/2}\big( p^\om_{\lfloor nt\rfloor }(0,g^\om_n(x))
 +  p^\om_{\lfloor nt\rfloor+1}(0,g^\om_n(x))\big) - 2 a^{-1} k^{(D)}_t(x)\Big| =0.
\end{equation}
For the continuous time random walk $Y$ we have
\begin{equation} \label{llt-ctsver}
  \lim_{n\to\infty}  
\sup_{x\in \bR^d} \sup_{t \ge T} 
 \Big|n^{d/2}q^\om_{nt}(0,g^\om_n(x))-  a^{-1} k^{(D)}_t(x)\Big| =0,
\end{equation}
where the constants $a$, $D$ are the same as for the myopic ant walk.
\end{theorem}

\ms We prove this theorem by establishing a parabolic Harnack
inequality (PHI) for solutions to the heat equation on $\Ci$.  (See
\cite{B1} for an elliptic Harnack inequality.)  This PHI implies
H\"older continuity of $p^\om_n(x,\cdot)$, and this enables us to
replace the weak convergence given by the CLT by pointwise
convergence.  In this paper we will concentrate on the proof of
\eqref{llt-pvi} -- the same arguments with only minor changes give
\eqref{llt-ctsver}.

Some of the results mentioned above, for random walks on percolation
clusters, have been extended to the `random conductance model', where
$\mu_{xy}$ are taken as i.i.d.r.v. in $[0,\infty)$ -- see \cite{BP, M,
SS}. In the case where the random conductors are bounded away from
zero and infinity, a local limit theorem follows by our methods -- see
Theorem \ref{thm:llt-rcm}.  If however the $\mu_{xy}$ have fat tails
at 0, then while a quenched invariance principle still holds, the
transition density does not have enough regularity for a local limit
theorem -- see Theorem 2.2 in \cite{BBHK}.

As an application of Theorem \ref{thm:llt-pi} we have the following
theorem on the Green's function $g_\om(x,y)$ on $\Ci$, defined 
(when $d \ge 3$) by
\begin{equation}\label{green-def}
 g_\om(x,y) = \int_0^\infty q^\om_t(x,y)dt. 
\end{equation}

\begin{theorem}\label{thm:green} Let $d \ge 3$.
(a) There exist constants $\delta, c_1, \dots c_4$, depending only on
$d$ and $p$, and r.v. $R_x$, $x \in \bZ^d$ 
such that
\begin{equation}\label{rxtail}
 \bP( R_x\ge n | x \in \Ci) \le c_1 e^{-{c_2} n^\delta},
 \end{equation}
for some $\delta=\delta(d,p)$, and non-random constants $c_i=c_i(d,p)$
such that
\begin{equation}\label{green-bds}
   \frac{c_3}{|x-y|^{d-2}} \le g_\om(x,y) \le \frac{c_4}{|x-y|^{d-2}} \qq 
\hbox{ if  } |x-y| \ge  R_x \wedge R_y.
\end{equation}
(b) There exists a constant $C=\Gamma(\frac{d}{2}-1)/(2\pi^{d/2}aD)>0$ such that 
for any $\eps>0$ there exists $M=M(\eps,\om)$ such that
on $\{ 0 \in \Ci\}$,
\begin{equation}\label{greenlim2}
 \frac{(1-\eps) C}{|x|^{d-2}} \le g_\om(0,x) \le 
 \frac{(1+\eps) C}{|x|^{d-2}} \q \hbox{ for } |x|>M(\om).
\end{equation}
(c) We have
\begin{equation}
 \lim_{|x| \to \infty} |x|^{2-d} \bE( g_\om(0,x) | 0 \in \Ci)
 = C.
 \end{equation}
\end{theorem}

\med {\bf Remark.} While \eqref{rxtail} gives
 good control of the tail of the random
variables $R_x$ in \eqref{green-bds}, we do not have any bounds 
on the tail of the r.v. $M$ in \eqref{greenlim2}. This is because
the proof of  \eqref{greenlim2} relies on the invariance principles in
\cite{SS, BB, MP}, and these do not give a rate of convergence. 

\ms

In Section 2 we indicate how the heat kernel estimates obtained in
\cite{B1} can be extended to discrete time, and also to variants of
the basic SRW $X$. In Section 3 we prove the PHI for $\Ci$ using the
`balayage' argument introduced in \cite{BBCK}. 
In the Appendix we give a self-contained proof of the key equation
in the simple fully discrete context of this section.
In Section 4 we show
that if the PHI and CLT hold for a suitably regular subgraph $\sG$ of
$\bZ^d$, then a local limit theorem holds.  In Section 5 we verify
these conditions for percolation, and prove Theorem \ref{thm:llt-pi}.
In Section 6, using the heat kernel bounds for $q^\om_t$ and the 
local limit theorem, we obtain Theorem \ref{thm:green}.

We write $c, c'$ for positive constants, which may change on
each appearance, and $c_i$ for constants which are fixed within
each argument. We occasionally use notation such as  $c_{1.2.1}$
to refer to constant $c_1$ in Theorem 1.2.

\section{ Discrete and continuous time walks}\label{sec2}

Let $\Gam=(G,E)$ be an infinite, connected graph with uniformly
bounded vertex degree. We write $d$ for the graph metric, 
and $B_d(x,r)=\{y : d(x,y) < r \}$ for balls with respect to $d$.
Given $A \subset G$, we write 
$\pd A$ for the external boundary of $A$
(so $y \in \pd A$ if and only if $y \in G-A$ and
there exists $x \in A$ with $x \sim y$.)
We set $\ol A = A \cup \pd A$.

Let $\mu_{xy}$ be `bond conductivities' on $\Gam$.
Thus $\mu_{xy}$ is defined for all $(x,y) \in G \times G$. 
We assume that $\mu_{xy}=\mu_{yx}$ for all
$x, y \in G$, and that $\mu_{xy}=0$ if $\{x,y\} \not\in E$ and
$x \neq y$. 
We assume that the conductivities on edges with distinct 
endpoints are bounded away from
$0$ and infinity, so that there exists a constant $C_M$ such that
\begin{equation}\label{mulbub}
  0< C_M^{-1} \le \mu_{xy} \le C_M  \qq \hbox{whenever }  x \sim y, \,
 x \neq y.
\end{equation}
We also assume that 
\begin{equation}\label{muxxub}
  0 \le \mu_{xx} \le C_M, \q  \hbox{for }  x\in G;
\end{equation}
we allow the possibility that $\mu_{xx}>0$ so as to be 
able to handle `blind ants' as in \cite{Ma}. 
We define $\mu_x=\mu(\{x\})=\sum_{y\in G} \mu_{xy}$, and extend $\mu$
to a measure on $G$. 
The pair $(\Gam,\mu)$ is often called a {\sl weighted graph}.
We assume that there exist $d\ge 1$ and $C_U$ such that
\begin{equation} \label{vub}
 \mu(B_d(x,r)) \le C_U r^d, \q  r\ge 1, \, x \in G.
\end{equation}
The standard discrete time SRW $X$ on $(\Gam,\mu)$ 
is the Markov chain $X=(X_n, n \in \bZ_+, P^x, x \in G)$ 
with transition probabilities $P_X(x,y)$ given by \eqref{PXdef}.
Since we allow $\mu_{xx}>0$, $X$ can 
jump from a vertex $x$ to itself.
We define the discrete time heat kernel on
$(\Gam,\mu)$ by
\begin{equation} 
 p_n(x,y) = \frac{P^x(X_n=y)}{\mu_x}. 
\end{equation}
 Let
\begin{equation} 
  \sL f (x) = \mu_x^{-1} \sum_y \mu_{xy} (f(y)-f(x)).
\end{equation}
One may also look at the continuous time SRW on 
 $(\Gam,\mu)$, which is the Markov process
$Y=(Y_t , t \in [0,\infty), \tilde P^x, x \in G)$,
with generator $\sL$. We define the 
(continuous time) heat kernel on $(\Gam, \mu)$ by
\begin{equation} 
 q_t(x,y) = \frac{\tilde P^x(Y_t=y)}{\mu_x}. 
\end{equation}
The continuous time heat kernel is a smoother object
that the discrete time one, and is often slightly 
simpler to handle.
Note that $p_n$ and $q_t $ satisfy
\begin{align*}
 p_{n+1}(x,y) -p_n(x,y) = \sL p_n(x,y), \qq
 \frac{\pd q_t(x,y)}{\pd t}  = \sL q_t(x,y).
\end{align*}
We remark that $Y$ can be constructed from $X$ by
making $Y$ follow the same trajectory as $X$, but
at times given by independent mean 1 exponential r.v.
More precisely, if $M_t$ is a rate 1 Poisson process, we set
$Y_t = X_{M_t}$, $t\ge 0$.
Define also the quadratic form
\begin{equation} 
  \sE(f,g) 
= \half \sum_x \sum_y \mu_{xy} (f(y)-f(x))(g(y)-g(x)).
\end{equation}

\cite{B1} studied the continuous time random walk $Y$ and the heat
kernel $q_t(x,y)$ on percolation clusters, in the case when
$\mu_{xy}=1$ whenever $\{x,y\}$ is an open edge, and $\mu_{xy}=0$
otherwise.  It was remarked in \cite{B1} that the same arguments work
for the discrete time heat kernel, but no details were given.  Since
some of the applications of \cite{B1} do use the discrete time
estimates, and as we shall also make use of these in this paper, we give
details of the changes needed to obtain these bounds.

In general terms, \cite{B1} uses two kinds of arguments to obtain the
bounds on $q_t(x,y)$.  One kind (see for example Lemma 3.5 or
Proposition 3.7) is probabilistic, and to adapt it to the discrete
time process $X$ requires very little work. The second kind uses
differential inequalities, and here one does have to be more careful,
since these usually have a more complicated form in discrete time.

We now recall some further definitions from \cite{B1}.

\med {\bf Definition} 
Let $C_V$, $C_P$, and $C_W\ge 1$ be fixed constants.
We say $B_d(x,r)$ is $(C_V,C_P,C_W)$--{\sl good} if:
\begin{equation}\label{vlb}
C_V r^d \le  \mu(B_d(x,r)),
\end{equation}
and the weak Poincar\'e inequality
\begin{equation}\label{wpi} 
\sum_{y \in B_d(x,r)} (f(y)-\ol{f}_{B_d(x,r)} )^2 \mu_y 
  \le C_P r^2 
\sum_{y,z  \in B_d(x, C_W r), z \sim y} |f(y)-f(z)|^2 \mu_{yz}
\end{equation}
holds for every $f:B_d(x,C_W r)\to \bR$. 
(Here $\ol f_{B_d(x,r)}$ is the value which minimises the
left hand side of (\ref{wpi})).

We say $B_d(x,R)$ is  $(C_V,C_P,C_W)$--{\sl very good} if there exists 
$N_B=N_{B_d(x,R)}\le R^{1/(d+2)}$
such that $B_d(y,r)$ is good whenever $B_d(y,r) \subseteq B_d(x,R)$, and 
$N_B \le r \le R$. We can always assume that $N_B \ge 1$. 
Usually the values of  $C_V,C_P,C_W$ will be clear from the context
and we will just use the terms `good' and `very good'.
(In fact the condition that $N_B \le  R^{1/(d+2)}$ is not used in this paper,
since whenever we use the condition `very good' we will impose 
a stronger condition on $N_B$).

\ms
 From now on in the section we fix 
$d\ge 2$, $C_M$, $C_V$, $C_P$, and $C_W$, and take
$(\Gam,\mu)=(G,E,\mu)$ to satisfy (\ref{vub}).
If $f(n,x)$ is a function on $\bZ_+ \times G$,
we write
\begin{equation}\label{fhat-def}
 \hat f(n,x) = f(n+1,x) + f(n,x),
\end{equation}
and in particular, 
to deal with the problem of bipartite graphs, 
we consider
\begin{equation} 
 \hat p_n(x,y) = p_{n+1}(x,y) + p_n(x,y).
\end{equation}

\smallskip
The following Theorem summarizes the bounds 
on $q$ and $p$ that will be used in the proof of the PHI and local
limit theorem. 

\begin{theorem}\label{basichk}
Assume that \eqref{mulbub}, \eqref{muxxub} and \eqref{vub} hold.
Let $x_0\in G$. 
Suppose that $R_1\ge 16$ and $B_d(x_0,R_1)$
is very good with $N_{B_d(x_0,R_1)}^{2d+4} \le R_1/(2 \log R_1)$.
Let $x_1 \in B_d(x_0, R_1/3)$.
Let $R \log R=R_1$, $T=R^2$, $B= B_d(x_1, R)$,
and $q^B_t(x,y)$, $p^B_n(x,y)$ be the heat kernels for the processes
$Y$ and $X$ killed on exiting from $B$.
Then
\begin{align}\label{bhk1}
 q^B_t(x,y) &\ge  c_1 T^{-d/2}, \q 
\hbox{ if } x, y \in B_d(x_1,3R/4), \q \fract14 T \le t \le T,\\
\label{bhk2}
 q_t(x,y) &\le  c_2 T^{-d/2}, \q 
\hbox{ if } x, y \in B_d(x_1,R), \q \fract14 T \le t \le T, \\
\label{bhk3}
 q_t(x,y) &\le  c_2 T^{-d/2}, \q 
\hbox{ if } x \in B_d(x_1,R/2),\q d(x,y) \ge R/8,
\q 0 \le t \le T,
\end{align}
and
\begin{align}\label{bhk4}
 p^B_{n+1}(x,y)+  &p^B_n(x,y)  \ge  c_1 T^{-d/2}, \q 
\hbox{ if } x, y \in B_d(x_1,3R/4), \q \fract14 T \le n \le T, \\
\label{bhk5}
 p_n(x,y) &\le  c_2 T^{-d/2}, \q 
\hbox{ if } x, y \in B_d(x_1,R), \q \fract14 T \le n \le T, \\
\label{bhk6}
 p_n(x,y) &\le  c_2 T^{-d/2}, \q 
\hbox{ if } x \in B_d(x_1,R/2),\q d(x,y) \ge R/8,
\q 0 \le n \le T.
\end{align}
\end{theorem}

\medskip
To prove this theorem we extend the bounds proved in \cite{B1}
for the continuous time simple random walk on $(\Gam,\mu)$
to the slightly more general random walks $X$ and $Y$ defined above.

\begin{theorem} \label{disc}
(a) Assume that \eqref{mulbub}, \eqref{muxxub} and \eqref{vub} hold.
Then the bounds in Proposition 3.1, Proposition 3.7,
Theorem 3.8, and Proposition 5.1-- Lemma 5.8 of \cite{B1}
all hold for $\hat p_n(x,y)$ as well as $q_t(x,y)$. \\
(b) In particular (see Theorem 5.7) let $x \in G$ and suppose that 
there exists $R_0=R_0(x)$ such that $B(x,R)$ is very good 
with $N_{B(x,R)}^{3(d+2)} \le R$ for each 
$R\ge R_0$. There exist constants $c_i$ 
such that if $n$ satisfies
$n \ge R_0^{2/3}$
then
\begin{equation}\label{ggub}
p_n(x,y) \le  c_1 n^{-d/2} e^{ -c_2 d(x,y)^2/n},  \qq
d(x,y) \le n, 
\end{equation}
and 
\begin{equation}\label{gglb}
p_n(x,y)+p_{n+1}(x,y) \ge  c_3 n^{-d/2} e^{ -c_4 d(x,y)^2/n},  \qq
 d(x,y)^{3/2} \le n. 
\end{equation}
(c) Similar bounds to those in \eqref{ggub}, \eqref{gglb} hold for 
$q_t(x,y)$.
\end{theorem}

\med{\bf Remark.} Note that we do not give in (b) Gaussian lower bounds
in the range $d(x,y) \le n < d(x,y)^{3/2}$.  However, as in
\cite[Theorem 5.7]{B1}, Gaussian lower bounds 
on $p_n$ and $q_t$ will hold in this range of values if a further
condition `exceedingly good' is imposed on $B(x,R)$ for all $R \ge
R_0$. We do not give further details here for two reasons; first the
 `exceedingly good' condition is rather complicated (see 
\cite[Definition 5.4]{B1}), and second the lower bounds in this range
have few applications.

\proof We only indicate the places where
changes in the arguments of \cite{B1} are needed.

First, let $\mu^0_{xy} =1$ if $\{x,y\} \in E$, and 0 otherwise. 
Then \eqref{mulbub} implies that 
if $\sE^0$ is the
quadratic form associated with $(\mu^0_{xy})$, then
\begin{equation} 
  c_1 \sE^0(f,f) \le \sE(f,f) \le  c_2 \sE^0(f,f)
\end{equation}
for all $f$ for which either expression is finite.  This means that
the weak Poincar\'e inequality for $\sE^0$ implies one (with a
different constant $C_P$) for $\sE$.  Using this, the arguments in
Section 3--5 of \cite{B1} go through essentially unchanged to give the
bounds for the continuous time heat kernel on $(\Gam,\mu)$.

More has to be said about the discrete time case. The 
argument in  \cite[Proposition 3.1]{B1} uses the equality
\begin{equation*} 
  \frac{\pd }{\pd t}q_{2t}(x_1,x_1) =- 2 \sE(q_t,q_t).  
\end{equation*}
Instead, in discrete time, we set
$f_n(x)= \hat p_n(x_1,x)$ and use the easily verified 
relation
\begin{equation} 
  \hat p_{2n+2}(x_1,x_1)- \hat p_{2n}(x_1,x_1)
 = -\sE(f_n, f_n).
\end{equation}
 Given this, the argument of  \cite[Proposition 3.1]{B1}
now goes through to give an upper bound on
$\hat p_n(x,x)$, and hence on $p_n(x,x)$.
A global upper bound, as in \cite[Corollary 3.2]{B1}, follows
since, taking $k$ to be an integer close to $n/2$,
\begin{align*} 
 p_n(x,y) = \sum_z p_k(x,z) p_{n-k}(y,z)
 &\le (\sum_z  p_k(x,z)^2)^{1/2} (\sum_z  p_{n-k}(y,z)^2)^{1/2} \\
 &= p_{2k}(x,x)^{1/2} p_{2n-2k}(y,y)^{1/2}.
\end{align*}

To obtain better bounds for $x,y$ far apart, \cite{B1}
used a method of Bass and Nash -- see \cite{Bas, Nash}.
This does not seem to transfer easily to discrete
time. For a process $Z$, write 
$\tau_Z(x,r)= \inf\{t: d(Z_t,x) \ge r\}$.
The key bound in continuous time is given
in  \cite[Lemma 3.5]{B1}, where it is proved that
if $B=B(x_0,R)$ is very good, then
\begin{equation}\label{l35}
 P^x( \tau_Y(x,r) \le t) \le \frac12 + \frac{ct}{r^2},
\hbox{ if } x \in B(x_0,2R/3), \q 0\le t \le c R^2/\log R,
\end{equation}
provided $c N_B^d (\log N_B)^{1/2} \le r \le R$.
(Here $N_B$ is the number given in the definition of `very good'.)
Recall that we can write $Y_t = X_{M_t}$, where $M$ is a rate 1 
Poisson process independent of $X$. So, 
$$ P^x( \tau_X(x,r) < t )P^x( M_{2t} > t)=
P^x( \tau_X(x,r) < t , M_{2t} >t)
 \le P(\tau_Y(x,r) < 2t). $$
Since $P( M_{2t} > t) \ge 3/4$ for $t\ge c$, we obtain
 \begin{equation} \label{xmeanub}
  P^x( \tau_X(x,r) < t ) \le \frac23 +  \frac{c't}{r^2}. 
\end{equation}
Using (\ref{xmeanub}) the remainder of the arguments
of Section  3 of \cite{B1} now follow through to give
the large deviation estimate Proposition 3.7 and 
the Gaussian upper bound Theorem 3.8.

The next use of differential inequalities in \cite{B1} is in
Proposition 5.1, where a technique of Fabes and Stroock \cite{FS} is
used. Let $B= B_d(x_1,R)$ be a ball in $G$, and $\vp: G \to \bR$, with
$\vp(x)>0$ for $x \in B$ and $\vp=0$ on $G-B$. Set
\begin{equation*} 
 V_0 = \sum_{x \in B} \vp(x) \mu_x.
\end{equation*}
Let $g_n(x)=\hat p_n(x_1,x)$, and 
\begin{equation} 
  H_n = V_0^{-1} \sum_{x \in B} \log( g_n(x)) \vp(x) \mu_x.
\end{equation}
We need to take $n \ge R$ here, so that $g_n(x)>0$
for all $x \in B$.
Using Jensen's inequality, and recalling that
$P_X(x,y)=\mu_{xy}/\mu_x$,
\begin{align}\label{discfs}
 \nonumber
 H_{n+1}-H_n &= \sum_{x \in B}  \log( g_{n+1}(x)/g_n(x)) \vp(x) \mu_x\\
  \nonumber
  &= \sum_{x \in G} \vp(x) \mu_x 
   \log \Big(\sum_{y\in G} P_X(x,y)  g_n(y)/g_n(x) \Big) \\
 \nonumber
   &\ge \sum_{x \in G} \vp(x) \mu_x 
    \sum_{y\in G} P_X(x,y) \log( g_n(y)/g_n(x)) \\
    \nonumber
  &=  \sum_{x \in G}  \sum_{y\in G}
  \vp(x) \mu_{xy} (\log g_n(y) -\log g_n(x)) \\
 &= -\half  \sum_{x \in G}  \sum_{y\in G} (\vp(y)-\vp(x))
 (\log g_n(y) -\log g_n(x)) \mu_{xy}.
\end{align}
Given (\ref{discfs}), the arguments on p.~3071-3073
of \cite{B1} give the basic 
`near diagonal' lower bound in \cite[Proposition 5.1]{B1},
for $\hat p_n(x,y)$.
The remainder of the arguments in Section 5 of \cite{B1}
can now be carried through. \qed

\sm {\it Proof of Theorem \ref{basichk}}.
This follows from Theorem \ref{disc}, using the fact that
Theorem 3.8 and Lemma 5.8 of \cite{B1} hold. \qed

\section { Parabolic Harnack Inequality} \label{sec-phi}

In this section we continue with the notation and hypotheses of
Section 2. Our first main result, Theorem \ref{hk2phi}, is a parabolic
Harnack inequality. Then, in Proposition \ref{prop:hc} we show that
solutions to the heat equation are H\"older continuous; this
result then provides the key to the local limit theorem proved
in the next section.

Let 
$$ Q(x,R,T)= (0,T] \times B_d(x,R), $$
and
$$ Q_-(x,R,T)=[\tfrac14 T, \tfrac12 T] \times B_d(x,\half R),
\q 
Q_+(x,R,T)=[\tfrac34 T,  T] \times B_d(x,\half R). $$
We use the notation $t+Q(x,R,T)= (t,t+T) \times B_d(x,R)$.
We say that a function $u(n,x)$ is {\sl caloric} on
$Q$ if $u$ is defined on 
$\ol Q = ([0,T]\cap \bZ) \times \ol B_d(x,R)$, and
\begin{equation}\label{calor}
  u(n+1,x) - u(n,x) = \sL u(n,x) \q \hbox { for } 
 0 \le n \le T-1, \, x \in B_d(x,R). 
\end{equation}
We say the parabolic Harnack inequality (PHI)
 holds with constant $C_H$ for
$Q=Q(x,R,T)$ if whenever $u=u(n,x)$ is non-negative and
caloric on $Q$, then
\begin{equation}\label{d-phi}
\sup_{(n,x) \in Q_-} \hat u(n,x) 
\le C_H \inf_{(n,x) \in Q_+} \hat u(n,x).
\end{equation}
The PHI in continuous time takes a similar form, except
that caloric functions satisfy
$$ \frac{\pd u}{\pd t} = \sL u, $$
and (\ref{d-phi}) is replaced by 
$\sup_{Q_-} u \le C_H \inf_{Q_+} u$.

\smallskip
We now show that the heat kernel bounds in Theorem \ref{basichk}
lead to a PHI.

\begin{theorem}\label{hk2phi} 
Let $x_0\in G$.  Suppose that $R_1\ge 16$ and $B_d(x_0,R_1)$ is 
$(C_V,C_P,C_W)$--very good
with $N_{B_d(x_0,R)}^{2d+4} \le R_1/(2 \log R_1)$.  
Let $x_1 \in B_d(x_0,R_1/3)$, and $R\log R=R_1$. Then there exists a 
constant $C_H$ such that the PHI (in both discrete and continuous time
settings) holds with constant $C_H$ for $Q(x_1,R,R^2)$.
\end{theorem}

\sm {\bf Remark.} 
The condition $R_1=R \log R$ here is not necessarily best possible.

\proof 
We use the balayage argument introduced in \cite{BBCK}
-- see also \cite{BBK2} for the argument in a graph setting. 
Let $T=R^2$, and write:
$$ B_0=B_d(x_1,R/2), \q B_1=B_d(x_1, 2R/3), \q B = B_d(x_1,R), $$
and
$$ Q=Q(x_1,R,T)=[0,T]\times B, \q E=(0,T]\times B_1. $$

We begin with the discrete time case.
Let $u(n,x)$ be non-negative and caloric on $Q$.
We consider the space-time process $Z$ on $\bZ \times G$
given by $Z_n=(I_n, X_n)$, where $X$ is the SRW on $\Gamma$,
$I_n=I_0-n$,
and $Z_0=(I_0,X_0)$ is the starting point of the space time process.
Define the r\'eduite $u_E$ by
\begin{equation} \label{red}
 u_E(n,x)= E^x \big( u(n-T_E, X_{T_E}); T_E< \tau_Q \big),   
\end{equation}
where $T_E$ is the hitting time of $E$ by $Z$, and $\tau_Q$
the exit time by $Z$ from $Q$.
So $u_E=u$ on $E$, $u_E=0$ on $Q^c$, and $u_E\le u$ on $Q-E$. 
As the process $Z$ has a dual, 
the balayage formula of Chapter VI of \cite{BG} holds
and we can write
\begin{equation}\label{bal1}
u_E(n,x)=  \int_E p^B_{n-r}(x,y) \nu_E(dr,dy), \q (n,x) \in Q, 
\end{equation}
for a suitable measure $\nu_E$. 
Here $p^B_n(x,y)$ is the transition density of the 
process $X$ killed on exiting from $B$. 

In this simple discrete setup we can write things more
explicitly. Set
\begin{equation}\label{J-def}
 J f(x) = 
\begin{cases}
 \sum_{y \in B} \frac{\mu_{xy}}{\mu_y} f(y), 
 &\hbox{ if } x \in B_1, \\
 0,  & \hbox{ if } x \in B-B_1. 
\end{cases}
\end{equation}
Then we have for $x \in B$,
\begin{equation}\label{bal-d1}
 u_E(n,x) =  \sum_{y \in B} p^B_n(x,y) u(0,y)\mu_y 
  +\sum_{y \in B} \sum_{r=2}^n p^B_{n-r}(x,y) k(r,y) \mu_y,
\end{equation}
where for $r \ge 2$
\begin{equation}\label{bal-k1}
   k(r,y) = J (u(r-1,\cdot)-u_E(r-1,\cdot))(y).
\end{equation}
See the appendix for a self-contained proof of \eqref{bal-d1}
and \eqref{bal-k1}.  

Since $u=u_E$ on $E$, if $r\ge 2$ then \eqref{bal-k1} implies
that $k(r,y)=0$ unless $y \in \pd (B-B_1)$.
Adding \eqref{bal-d1} for $u(n,x)$ and $u(n+1,x)$, and using the 
fact that  $k(n+1,x)=0$ for $x \in B_0$,
we obtain, for $x \in B_0$,
\begin{equation}\label{bal-d2}
 \hat u_E(n,x) = 
\sum_{y \in B_1} \sum_{r=1}^n \hat p^B_{n-r}(x,y) k(r,y) \mu_y.
\end{equation}

Now let $(n_1,y_1) \in Q_-$ and $(n_2,y_2) \in Q_+$. 
Since $(n_i,y_i)\in E$ for $i=1,2$, 
we have $u_E(n_i,y_i)=u(n_i,y_i)$, and so \eqref{bal-d2} holds.
By Theorem \ref{basichk} we have, writing
$A=  \pd (B-B_1)$,
\begin{align*}
 \hat p^B_{n_2-r}(x,y) & \ge  c_1 T^{-d/2}
  \q \hbox{ for } x, y \in B_1, \, 0\le r \le T/2, \\
 \hat p_{r}(x,y) &\le  c_2  T^{-d/2}
  \q \hbox{ for } x, y \in B_1, \, T/4 \le r \le T/2, \\
  \hat p_{n_1-r}(x,y) &\le c_2 T^{-d/2}
  \q \hbox{ for } x \in B_0, \, y \in A, \, 0 < r \le n_1.
\end{align*}
Substituting these bounds in \eqref{bal-d2}, 

\begin{align*}
  \hat u(n_2,y_2) &=
 \sum_{y \in B_1} \hat p^B_{n_2}(y_2,y) u(0,y)\mu_y
+  \sum_{y \in A} \sum_{r=2}^{n_2} \hat p^B_{n_2-s}(y_2,y)k(r,y) \mu_y  \\
 &\ge 
 \sum_{y \in B_1} \hat p^B_{n_2}(y_2,y) u(0,y)\mu_y
+  \sum_{y \in A} \sum_{r=2}^{n_1} \hat p^B_{n_2-s}(y_2,y)k(r,y) \mu_y  \\
 &\ge 
 \sum_{y \in B_1} c_1 T^{-d/2} u(0,y)\mu_y
+  \sum_{y \in A} \sum_{r=2}^{n_1} c_1 T^{-d/2} k(r,y) \mu_y  \\
 &\ge 
 \sum_{y \in B_1} c_1 c_2^{-1} \hat p^B_{n_1}(y_1,y) u(0,y)\mu_y
+  \sum_{y \in A} \sum_{r=2}^{n_1}c_1 c_2^{-1} \hat p^B_{n_1-s}(y_1,y)k(r,y) 
 \mu_y \\
 &=  c_1 c_2^{-1} \hat u(n_1,y_1),
\end{align*}
 which proves the PHI. 

The proof is similar in the continuous time case. 
The balayage formula takes the form 
\begin{equation}\label{bal-cts}
 u_E(t,x) = \sum_{y \in B} q^B_t(x,y) u(0,y)\mu_y
+   \sum_{y \in B_1}\int_0^t q^B_{t-s}(x,y)k(s,y) \mu_y ds,
\end{equation}
where $k(s,y)$ is zero if $y \in B-B_1$ and  
\begin{equation}\label{k-cts}
   k(s,y) = J (u(s,\cdot)-u_E(s,\cdot))(y), \q y \in B_1.
\end{equation}
(See \cite[Proposition 3.3]{BBK2}).
Using the bounds on $q^B_t$ in Theorem \ref{basichk}
then gives the PHI.
\qed

\sm {\bf Remark.} In \cite{B1} an elliptic Harnack inequality (EHI)
was proved for random walks on percolation clusters -- see Theorem
5.11. Since the PHI immediately implies the EHI, the argument above
gives an alternative, and simpler, proof of this result.

\ms It is well known that the PHI implies H\"older continuity of
caloric functions -- see for example Theorem 5.4.7 of \cite{SC}. But
since in our context the PHI does not hold for all balls, we give the
details of the proof.  In the next section we will just use this
result when the caloric function $u$ is either $q_t(x,y)$ or 
$\hat p_n(x,y)$.

\begin{proposition}\label{prop:hc}
Let $x_0 \in G$. Suppose that there exists $s(x_0)\ge 0$ such
that the PHI (with constant $C_H$)
holds for $Q(x_0,R,R^2)$ for $R \ge s(x_0)$. 
Let $\theta = \log( 2C_H/(2C_H-1))/\log 2$, and 
\begin{equation}
 \rho(x_0,x,y) = s(x_0) \vee d(x_0,x) \vee d(x_0,y).
\end{equation}
Let $r_0 \ge s(x_0)$, $t_0=r_0^2$, and suppose that 
$u=u(n,x)$ is caloric in $Q=Q(x_0,r_0,r_0^2)$.
Let $x_1,x_2\in B_d(x_0, \half r_0)$, and
$t_0-  \rho(x_0,x_1,x_2)^2 \le n_1, n_2 \le t_0-1$.
Then
\begin{equation}\label{u-hold}
 |\hat u(n_1,x_1)- \hat u(n_2,x_2)| 
\le c \Big(\frac{ \rho(x_0,x_1,x_2)}{t_0^{1/2}}\Big)^\theta
  \sup_{Q_+} |\hat u|.
\end{equation}
\end{proposition}

\proof We just give the discrete time argument --
the continuous time one is almost identical.
Set $r_k = 2^{-k} r_0$, and let
$$ Q(k)=  (t_0-r_k^2)+ Q(x_0,r_k,r_k^2). $$
Thus $Q_+(k)=Q(k+1)$.
Let $k$ be such that $r_k \ge s(x_0)$. Let
$\hat v$ be $\hat u$ normalised in $Q(k)$ so 
that $0\le \hat v \le 1$, and $\Osc(\hat v,Q(k))=1$.
(Here $\Osc(u,A) = \sup_Q u -\inf_A u$ is the oscillation
of $u$ on $A$).
Replacing $\hat v$ by $1-\hat v$ if necessary we can assume
$ \sup_{Q_-(k)} \hat v \ge \half$.
By the PHI,
$$ \half \le \sup_{Q_-(k)} \hat v \le C_H \inf_{Q_+(k)} \hat v, $$
and it follows that, if $\delta = (2 C_H)^{-1}$, then
\begin{equation} \label{oscin}
\Osc(\hat u,Q_+(k)) \le (1-\delta) \Osc(\hat u,Q(k)). 
\end{equation}

Now choose $m$ as large as possible so that
$r_m \ge  \rho(x_0,x,y)$.
Then applying (\ref{oscin}) in the chain of boxes
$Q(1) \supset Q(2) \supset \dots Q(m)$, we deduce that,
since $(x_i,n_i) \in Q(m)$, 
\begin{equation}\label{u-ba}
 |\hat u(n_1,x_1)- \hat u(n_2,x_2)| \le \Osc(\hat u,Q_m)
 \le (1-\delta)^{m-1} \Osc(\hat u,Q(1)). 
\end{equation}
Since $(1-\delta)^m \le c (r_0/t_0^{1/2})^\theta$, 
(\ref{u-hold}) follows from  (\ref{u-ba})  . \qed

\section{ Local limit theorem}

\medskip Now let $\sG \subset \bZ^d$, and let $d$ denote graph
distance in $\sG$, regarded as a subgraph of $\bZ^d$.
We assume $\sG$ is infinite and connected, and $0 \in \sG$. 
We define $\mu_{xy}$ as in Section 2 so that
\eqref{mulbub}, \eqref{muxxub} and \eqref{vub} hold, and
write $X=(X_n, n \in \bZ_+, P^x, x \in \sG)$  
for the associated simple random walk on $(\sG,\mu)$.
We write $|\cdot|_p$ for the $L^p$ norm in $\bR^d$;
$|\cdot|$ is the usual ($p=2$) Euclidean distance.

Recall that $k_t^{(D)}(x)$ is the Gaussian heat kernel in $\bR^d$ with
diffusion constant $D>0$ and 
let $X^{(n)}_t =  n^{-1/2} X_{\lfloor nt \rfloor}$.
For $x\in\bR^d$, set 
\begin{equation}
 H(x,r) = x + [-r,r]^d, \qq \Lam(x,r) = H(x,r) \cap \sG.
\end{equation}
In general $\Lam(x,r)$ will not be connected.
Let 
$$  \Lam_n(x,r) = \Lam(xn^{1/2},rn^{1/2}). $$
Choose a function $g_n: \bR^d \to \sG$ so that 
$g_n(x)$ is a closest point in $\sG$ to $n^{1/2} x$, 
in the  $|\cdot|_\infty$ norm.
(We can define $g_n$ by using some fixed ordering of
$\bZ^d$ to break ties.)

We now make the following assumption on the graph $\sG$
and the SRW $X$ on $\sG$.
Let  $x \in \bR^d$.

\begin{assumption}\label{ass1}
There exists
a constant $\delta>0$, and positive constants $D, C_H$, 
$C_i, a_\sG$ such that the following hold. 
\newline (a) (CLT for $X$). 
For any $y \in \bR^d$, $r>0$,
\begin{equation}
 P^0( X^{(n)}_{t} \in H(y,r) ) \to \int_{H(y,r)} k^{(D)}_t(y') dy'. 
\end{equation}
(b) There is a global upper heat kernel bound of the form
\[ p_k(0,y) \le C_2 k^{-d/2}, \q \hbox{ for all } y \in \sG, k \ge C_3. \]
(c) For each $y \in \sG$ there exists $s(y)<\infty$ such that
the PHI \eqref{d-phi} 
holds with constant $C_H$ for $Q(y,R,R^2)$ for $R\geq s(y)$. \\
(d) For any $r>0$
 \begin{equation}\label{a-eqn}
 \frac{\mu(\Lam_n(x,r))}{ (2 n^{1/2}r)^d} \to a_\sG \qq
 \hbox{ as } n \to \infty. 
\end{equation}
(e) For each $r>0$ there exists $n_0$ such
that, for $n\ge n_0$,
\[ |x'-y'|_\infty \le d(x',y') \le (C_1|x'-y'|_\infty) \vee n^{1/2-\delta}, \;\; 
\hbox{ for all } x',y' \in \Lam_n(x,r). \]
(f) $n^{-1/2} s(g_n(x)) \to 0$ as $n\to\infty$.
\end{assumption}

We remark that for any $x$ all these hold for $\bZ^d$: for 
the PHI see \cite{Del}.
We also remark that these assumptions are not independent; for example
the PHI in (c) implies an upper bound as in (b).
For the region  $Q(y,R,R^2)$ in (c) the space ball is in the graph
metric on $\sG$.

\smallskip
We write, for $t \in [0,\infty)$,
$$  \hat p_t(x,y) = \hat p_{\lfloor t \rfloor}(x,y)
  =  p_{\lfloor t \rfloor}(x,y)+ p_{\lfloor t \rfloor+1}(x,y). $$

\begin{theorem}\label{thm:llt}
Let $x \in \bR^d$ and $t>0$. Suppose Assumption \ref{ass1}
holds. Then
\begin{equation}
\label{llt1}
  \lim_{n \to \infty} n^{d/2} \hat p_{nt}(0,g_n(x)) 
=  2 a_\sG ^{-1} k^{(D)}_t(x).
\end{equation}
\end{theorem}

\proof Write $k_t$ for $k_t^{(D)}$.
Let $\theta$ be chosen as in Proposition \ref{prop:hc}.
Let $\eps \in (0, \half)$.  
Choose $\kappa>0$ such that
$ ( \kappa^\theta + \kappa)< \eps$. 
Write $\Lam_n=\Lam_n(x, \kappa)= \Lam(n^{1/2}x, n^{1/2}\kappa)$. 
Set
\begin{equation}
 J(n) =  P^0 \Big( n^{-1/2} X_{\lfloor nt \rfloor}  \in \Lam(x,\kappa) \Big) 
     + P^0\Big( n^{-1/2} X_{\lfloor nt \rfloor +1}  \in \Lam(x,\kappa) \Big) 
  - 2 \int_{\Lam(x,\kappa)} k_t(y) dy . \\
\end{equation}
Then 
\begin{align}
 \nonumber
J(n) &= \sum_{z \in \Lam_n} 
 \big(\hat p_{nt}(0,z) - \hat p_{nt}(0,g_n(x)) \big)\mu_z \\
  \label{j2}
 &\qq +  \mu(\Lam_n)  \hat p_{nt}(0,g_n(x)) -
  \mu(\Lam_n)n^{-d/2} a_\sG^{-1} 2 k_t(x)  \\
 \label{j3}
 &\qq \q 
 +  2 k_t(x) \big( \mu(\Lam_n)n^{-d/2} a_\sG^{-1} -2^d \kappa^d \big)\\
 \label{j4}
 &\qq \qq +   2 \int_{H(x,\kappa)} (k_t(x)-k_t(y))dy \\
 &= J_1(n) + J_2(n) +J_3(n) +J_4(n).
\nonumber
\end{align}

We now control the terms $J(n)$, $J_1(n)$, $J_3(n)$ and
$J_4(n)$. By Assumption \ref{ass1}
we can choose $n_1$ with $n_1^{-\delta} < 2 C_1 \kappa$ 
such that, for $ n \ge n_1$,
\begin{align}
 \label{jbnd}
 |J(n)| &\le \kappa^d \eps, \\
 \label{mubnd}
  \left| \frac{\mu(\Lam_n)}{a_\sG (2 n^{1/2}\kappa)^d}  - 1 \right|
  &\le \eps < \half, \\
  \label{pbnd}
 \sup_{k \ge \half nt, z \in \sG} \hat p_{k}(0,z)
  &\le c_1 (nt)^{-d/2}, \\
  \label{sgbnd}
  s(g_n(x)) n^{-1/2} &\le 2 C_1 \kappa.
\end{align}

We bound $J_1(n)$ by using the H\"older continuity of $\hat p$, 
which comes from the PHI and Proposition \ref{prop:hc}.
We begin by comparing $\Lam_n$ with balls in the $d$-metric.
Let $n \ge n_1$.
By \eqref{mubnd} $\mu(\Lam_n)>0$, so $g_n(x) \in \Lam_n$.
By Assumption  \ref{ass1}(e) there exists $n_2 \ge n_1$ such that,
if $n \ge n_2$ and $y \in \Lam_n$ then
\begin{equation*}
 d(y,g_n(x)) \le (C_1 |y-g_n(x)|_\infty) \vee n^{1/2-\delta} 
\le  n^{1/2}\big( (2C_1 \kappa) \vee n^{-\delta}\big)
 \le 2C_1 \kappa n^{1/2}.
\end{equation*}
So, writing  $B= B_d(g_n(x),2C_1 \kappa n^{1/2})$, 
$\Lam_n \subset B$ when $n \ge n_2$.
Thus we have, using (\ref{mubnd}),
\begin{align}\nonumber
 |J_1(n)| &\le \mu(\Lam_n) 
 \max_{z \in \Lam_n} |\hat p_{nt}(0,z) - \hat p_{nt}(0,g_n(x))|\\
 \label{j1ba}
 &\le  2 a_\sG  (2 n^{1/2}\kappa)^d 
 \max_{z \in B} |\hat p_{nt}(0,z) -  \hat p_{nt}(0,g_n(x))|. 
\end{align}
Using Assumption \ref{ass1}(c), Proposition \ref{prop:hc}
and then (\ref{pbnd}) and  (\ref{sgbnd}),
\begin{align}\nonumber
 \max_{z \in B} 
|\hat p_{nt}(0,z) -  \hat p_{nt}(0,g_n(x))| 
 &\le c \Big( \frac{s(g_n(x))\vee 2 C_1\kappa n^{1/2}}
      {(nt)^{1/2}} \Big)^\theta 
\sup_{k \ge \half nt, z \in \sG} \hat p_{k}(0,z) \\ 
\nonumber
 &\le c (nt)^{-d/2} 
 \Big( \frac{s(g_n(x)) n^{-1/2} \vee 2 C_1\kappa }
 {t^{1/2}} \Big)^\theta \\
\label{maxpbd}
  &\le c_2 t^{-(d+\theta)/2}  n^{-d/2} \kappa^{\theta}.
 \end{align}
Hence combining \eqref{j1ba} and \eqref{maxpbd}
\begin{equation}\label{jibz}
 |J_1(n)| \le  c_3 t^{-(d+\theta)/2} \kappa^{d+\theta}.
\end{equation}

We now control the other terms. Since $|\nabla k_t(x)| \le c_4 t^{-(d+1)/2}$,
\begin{equation}
 |J_4(n)| \le 2 |\Lam(x,\kappa)| c_4(t) (2\kappa) = \kappa^{d+1} c_5(t).
\end{equation}
For $J_3(n)$, using (\ref{mubnd}) and  (\ref{pbnd}), if
$n \ge n_2$ then
\begin{align*}
 J_3(n) &= 2  k_t(x)
 \big| \mu(\Lam_n)n^{-d/2} a_\sG^{-1} - 2^d \kappa^d \big|\\
  &= 2 k_t(x) 2^d \kappa^d  
 \Big| \frac{\mu(\Lam_n)}{a_\sG (2 n^{1/2}\kappa)^d} -1 \Big| 
 \le c_6(t)\kappa^d \eps. 
\end{align*}

Now write $\wt p_n = n^{d/2} \hat p_{nt}(0,g_n(x))$. 
Then for $n\geq n_2$
\begin{align*}
 |J_2(n)| 
 &= \mu(\Lam_n)|\hat p_{nt}(0,g_n(x)) - n^{-d/2} a_\sG^{-1} 2 k_t(x)| \\
 &= \frac{ \mu(\Lam_n)}{(2 n^{1/2} \kappa)^d}
 (2\kappa)^d   |\wt p_n  - 2 a_\sG^{-1} k_t(x)| 
 \ge \half a_\sG (2\kappa)^d  |\wt p_n  - 2 a_\sG^{-1} k_t(x)|.
 \end{align*}
So, 
\begin{align*}
 \half a_\sG (2\kappa)^d  | \wt p_n  - 2 a_\sG^{-1} k_t(x)|
 &\le | J(n)|+ |J_1(n)|+|J_3(n)|+|J_4(n)| \\
 &\le \kappa^d \eps + c_3 t^{-(d+\theta)/2} \kappa^{d+\theta}
  +   c_6(t)\kappa^{d} \eps  +  c_5(t) \kappa^{d+1} \\
 &\le c_7(t) \kappa^d ( \eps+ \kappa^\theta + \kappa) 
  \le 2 c_7(t) \kappa^d \eps.
\end{align*}
Thus for $n \ge n_2$,
\begin{equation}\label{ptil-b}
   | \wt p_n  - 2 a_\sG ^{-1} k_t(x)|
 \le c_8(t) \eps,
\end{equation}
which completes the proof. \qed

\begin{corollary}\label{u-t} 
Let  $0<T_1<T_2<\infty$. 
Suppose  Assumption \ref{ass1} holds, and in addition
that for each $H(y,r)$ the CLT in Assumption \ref{ass1}(a) holds
uniformly for $t \in [T_1, T_2]$.
Then 
\begin{equation}
\label{llt2}
  \lim_{n \to \infty} \sup_{T_1 \le t \le T_2}
 | n^{d/2} \hat p_{nt}(0,g_n(x)) - 2 a_\sG ^{-1} k^{(D)}_t(x)|=0.
\end{equation}
\end{corollary}

\proof The argument is the same as for the Theorem; all we
need do is to note that the constant $c_8(t)$ in
(\ref{ptil-b}) can be chosen to be bounded on $ [T_1, T_2]$.
\qed

If we slightly strengthen our assumptions, then 
we can obtain a uniform result in $x$. 

\begin{assumption}\label{assu}
(a) For any compact $I \subset (0,\infty)$, 
the CLT in  Assumption \ref{ass1}(a) holds uniformly
for $t \in I$. \\
(b) There exist $C_i$ such that
\begin{equation} \label{fgub}
\hat p_k(0,x) \le C_2 k^{-d/2} \exp( - C_4 d(0,x)^2 /k), \qq
\hbox{ for $k \ge C_3$ and $x \in \sG$}. 
\end{equation}
(c) Assumption \ref{ass1}(c) holds. \\
(d) Let $h(r)$ be the size of the biggest `hole'
in $\Lam(0,r)$. More precisely, 
$h(r)$ is the suprema of the
$r'$ such that $\Lam(y,r')= \emptyset$ for some $y \in H(0,r)$.
Then $\lim_{r \to \infty} h(r)/r= 0$. \\
(e) There exist constants $\delta$, $C_1$, $C_H$ such that
for each $x \in \bQ^d$ Assumption \ref{ass1}(d), (e) and (f) hold. 
\end{assumption}
Note that in discrete time we have $p_k(0,x)=0$ if
$d(0,x)>k$, so it is not necessary in \eqref{fgub} to consider separately
the case when $d(0,x) \gg k$.

\begin{theorem}\label{thm:unifllt}
Let $T_1>0$. Suppose Assumption \ref{assu} holds. Then
\begin{equation}
\label{ullt1}
  \lim_{n\to\infty}  
\sup_{x\in \bR^d} \sup_{t \ge T_1} 
 |n^{d/2} \hat p_{nt}(0,g_n(x)) - 2 a_\sG ^{-1} k_t^{(D)}(x)| =0.
\end{equation}
\end{theorem}

\proof As before we write $k_t=k_t^{(D)}$. Set
$$ w(n,t,x)=  |n^{d/2} \hat p_{nt}(0,g_n(x)) - 2 a_\sG ^{-1} k_t(x)|.$$
Let $\eps \in (0, \half)$. 
We begin by restricting to a compact set of $x$ and $t$.
Choose $n_1$ so that $n_1 T_1\ge C_3$, and $T_2> 1+ T_1$ such that
\begin{align*}
  2 a_\sG^{-1} k_{T_2}(0) +  C_2 T_2^{-d/2}  \le \eps.
\end{align*}
If $t\ge T_2$ then using Assumption \ref{ass1}(b), for $n \ge n_1$,
\begin{align*}
  w(n,t,x) 
 &\le n^{d/2} \hat p_{nt}(0,g_n(x)) +  2 a_\sG ^{-1} k_t(x) 
 \le  n^{d/2} C_2 (nt)^{-d/2} + 2 a_\sG ^{-1} k_t(0) \le \eps. 
\end{align*}
So we can restrict to $t \in [T_1,T_2]$.

Now choose $R>0$ so that  $h(r) \le \half r$ for $r \ge R$.
Let $|x| \ge R$ and  $t \in [T_1,T_2]$. Then
\begin{equation}\label{r-0}
2 a_\sG ^{-1} k_t(x) \le c T_1^{-d/2}  \exp( - R^2/2 T_2).
\end{equation}
We have 
$ |n^{1/2} x - g_n(x)|_\infty \le h(|x| n^{1/2}) \le \half |x| n^{1/2}$,
as $|x|n^{1/2}>R$ for all $n\geq 1$, and hence
$$ d(0, g_n(x)) \ge |g_n(x)|_\infty  \ge \half |x| n^{1/2}.$$
The Gaussian upper bound \eqref{fgub} yields 
\begin{equation}\label{pu-1}
 n^{d/2} \hat p_{nt}(0,g_n(x))
 \le c t^{-d/2} \exp (-c' |x|^2/t) 
\le c T_1^{-d/2} \exp( - c' R^2/T_2).
\end{equation}
We can choose $R$ large enough so the terms in
(\ref{r-0}) and (\ref{pu-1}) are smaller than
$\eps$. Thus $w(n,t,x)< \eps$ whenever $t>T_2$ or
$|x|>R$, and 
$n \ge n_1$.
Thus it remains to show that there exists $n_2$
such that for $n \ge n_2$,
$$ \sup_{ |x|\le R, T_1 \le t \le T_2} w(n,t,x) < \eps . $$

Now let $\kappa$ be chosen as in the proof of 
Theorem \ref{thm:llt}, and also such that
\begin{equation}\label{kap-cho} 
c_1 T_1^{-(d+\theta)/2} \kappa^{\theta} < \eps,
\end{equation}
where $c_1$ is the constant $c_3$ in \eqref{jibz}.
Let 
$\eta \in (0, \kappa) \cap  \bQ$. Set
$\sY = \{y\in \eta\bZ^d \cap B_R(0)\}$,
where $B_R(0)$ is the Euclidean ball centre 0 and radius $R$.
By  Theorem \ref{thm:llt} and Corollary \ref{u-t}
for each $y \in \sY$ there exists $n'_3(y)$ such that
\begin{equation}\label{wbnd-1}
   \sup_{T_1\le t\le T_2} w(n,t,y) \le \eps
\q \hbox{ for } n \ge n'_3(y).
\end{equation}
We can assume in addition that $n'_3(y)$ is greater than the $n_2=n_2(y)$
given by the proof of Theorem \ref{thm:llt}.
Let $n_4 = \max_{y \in \sY} n'_3(y)$.
Now let $x \in B_R(0)$, and write $y(x)$ for a closest
point (in the $|\cdot|_\infty$ norm) in $\sY$ to $x$: 
thus $|x - y(x)|_\infty \le \eta$.
Let $n \ge n_4 $. We have 
\begin{align}\label{app-1}
  |n^{d/2}\hat p_{nt}(0,g_n(x)) - 2 a_\sG^{-1}k_t(x)| &\le
  |n^{d/2}\hat p_{nt}(0,g_n(x)) - n^{d/2}\hat p_{nt}(0,g_n(y(x)))| \\
\label{app-2}
 & \qquad +  |n^{d/2}\hat p_{nt}(0,g_n(y(x)))- 2 a_\sG^{-1}k_t(y(x))| \\ 
 \label{app-3}
& \qquad \qquad + |2 a_\sG^{-1} k_t(y(x))- 2 a_\sG^{-1}k_t(x)|,
\end{align}
and it remains to bound the three terms (\ref{app-1}), 
(\ref{app-2}), (\ref{app-3}), which we denote $L_1, L_2, L_3$ respectively.
Since $\eta< \kappa$ and $n \ge n_4 \ge n_3(y(x))$,
we have the same bound for $L_1$ as in \eqref{maxpbd}, and obtain
\begin{align} \label{uub-c}
  L_1= |n^{d/2}\hat p_{nt}(0,g_n(x)) - n^{d/2} \hat p_{nt}(0,g_n(y(x)))|
&\le c_1 t^{-(d+\theta)/2} \eta^{\theta} \\
&\le  c_1 T_1^{-(d+\theta)/2} \eta^{\theta} < \eps ,
\end{align}
by \eqref{kap-cho}.
As $n \ge n_4$ and $y(x)\in \sY$, by (\ref{wbnd-1}) $L_2 < \eps$.
Finally,
$$ L_3 = |k_t(x)-k_t(y(x))| \le \eta d^{1/2} ||\grad k_t||_\infty
  \le c  \eta T_1^{-(d+1)/2}, $$
and choosing $\eta$ small enough this is less than $\eps$.
Thus we have $w(n,t,x)< 3 \eps$
for any $x \in B_R(0)$, $ t \in [T_1, T_2]$ and $n \ge n_4$,
completing the proof of the theorem. \qed

\medskip
In continuous time we replace $X$ by $Y$, $p_k(0,y)$ by $q_t(0,y)$,
and modify Assumptions \ref{ass1} and \ref{assu} accordingly. 
That is, in both Assumptions we replace the CLT for $X$
in (a) by a CLT for $Y$, replace $p_n$ in (b) by $q_t$, 
and require the continuous time version of the PHI in (c).   
The same
arguments then give a local limit theorem as follows.

\begin{theorem}\label{thm:unifllt-cts}
Let $T_1>0$. Suppose Assumption \ref{assu} 
(modified as above for the continuous time case) holds. Then
\begin{equation}
\label{ullt1-cts}
  \lim_{n\to\infty}  
\sup_{x\in \bR^d} \sup_{t \ge T_1} 
 |n^{d/2} q_{nt}(0,g_n(x)) -  a_\sG ^{-1} k_t^{(D)}(x)| =0.
\end{equation}
\end{theorem}

\section{Application to percolation clusters}\label{sec-perc}

We now let $(\Omega, \bP)$ be a probability space
carrying a supercritical bond percolation process on $\bZ^d$.
As in the Introduction
we write $\Ci=\Ci(\om)$ for the infinite cluster. 
Let $\bP_0(\cdot) = \bP(\cdot|  0 \in \Ci)$.
Let $x \sim y$. We set $\mu_{xy}(\om)=1$ if the edge $\{x,y\}$ is open
and  $\mu_{xy}(\om)=0$ otherwise.
In the physics literature one finds two common choices
of random walks on $\Ci$, called the `myopic ant' and 'blind ant' walks,
which we denote $X^M$ and $X^B$ respectively.
For the myopic walk we set
\begin{align*}
  \mu^M_{xy} &= \mu_{xy}, \q y \neq x, \\
  \mu^M_{xx} &= 0,
\end{align*}
and for each $\om \in \Omega$ we then take 
$X^M=(X^M_n, n \in \bZ_+, P^x_\om, x \in \Ci(\om))$ to be the random 
walk on the graph $(\Ci(\om), \mu^M(\om))$.
Thus $X^M$ jumps with equal probability from $x$ 
along any of the open bonds adjacent to $x$. 
The second choice (`the blind ant') 
is to take 
\begin{align*}
  \mu^B_{xy} &= \mu_{xy}, \q y \neq x, \\
  \mu^B_{xx} &= 2d- \mu_x,
\end{align*}
and take $X^B$  to be the random 
walk on the graph $(\Ci(\om), \mu^B(\om))$.
This walk attempts to jump with probability $1/2d$ in 
each direction, but the jump is suppressed if the bond is not open.  
By Theorem \ref{disc} the same transition density 
bounds hold for these two processes. 
Since these two processes are time changes of each other, 
an invariance principle for one quickly leads to one for the
other -- see for example \cite[Lemma 6.4]{BB}.
 
In what follows we take $X$ to be either of the two walks given
above. We write $p^\om_n(x,y)$ for its transition density, and as
before we set $\hat p^\om_n(x,y) = p^\om_n(x,y) + p^\om_{n+1}(x,y)$.
We begin by summarizing the heat kernel bounds on $p^\om_n(x,y)$. 

\begin{theorem}\label{disc-hkb-p}
There exists $\eta=\eta(d)>0$ and  constants $c_i=c_i(d,p)$ and 
r.v. $V_x, x\in \bZ^d$, such that
\begin{equation}\label{vxub}
\bP( V_x(\omega) \ge n) \le c \exp(-c n^{\eta}), 
\end{equation}
and if $ n \ge c |x-y| \vee V_x$ then
\begin{align}\label{pn-gb}
 c_1 n^{-d/2} e^{-c_2 |x-y|^2/n} \le 
\hat p^\om_n(x,y) \le c_3 n^{-d/2} e^{-c_4 |x-y|^2/n}.
\end{align}
Further if $n \ge c |x-y|$ then
\begin{align}\label{Epn-gb}
 c_1 n^{-d/2} e^{-c_2 |x-y|^2/n} \le 
\bE( \hat p^\om_n(x,y)| x,y \in \Ci) \le c_3 n^{-d/2} e^{-c_4 |x-y|^2/n}.
\end{align}
\end{theorem}

\proof 
This follows from Theorem \ref{disc}(a), and the arguments
in \cite{B1}, Section 6.  \qed

\sms We now give the local limit theorem. 
As in Section 4 we write $g^\om_n(x)$ for a closest point
in $\Ci$ to $n^{1/2} x$, set
$\Lam(x,r)=\Lam(x,r)(\om) = \Ci(\om) \cap H(x,r)$, 
and write $h_\om(r)$ for the largest hole in $\Lam(0,r)$.

\begin{theorem}\label{thm:llt-p}
Let $T_1>0$. Then there exist constants $a$, $D$ such that
$\bP_0$-a.s.,
\begin{equation}
\label{llt-pv}
  \lim_{n\to\infty}  
\sup_{x\in \bR^d} \sup_{t \ge T_1} 
 |n^{d/2} \hat{p}^\om_{nt}(0,g^\om_n(x)) - 2 a^{-1} k^{(D)}_t(x)| =0.
\end{equation}
\end{theorem}

\noindent In view of Theorem \ref{thm:unifllt} it is enough to prove
that, $\bP_0$-a.s., the cluster $\Ci(\om)$ and process $X$
satisfy Assumption
\ref{assu}. Note that since we apply Theorem \ref{thm:unifllt}
separately to each graph $\Ci(\om)$, it is not necessary that the
constants $C_i$ in Assumption \ref{assu} should be uniform in $\om$
-- in fact, it is clear that the constant $C_3$ in \eqref{fgub} 
cannot be taken independent of $\om$.

\begin{lemma}\label{easy-a}
 (a) There exist constants $\delta$, $C_\cdot$
such that Assumption \ref{assu} (a), (b), (c)
all hold $\bP_0$-a.s. \\
(b) Let $x \in \bR^d$. Then Assumption \ref{ass1}(e)
holds  $\bP_0$-a.s.
\end{lemma}

\proof (a) The CLT holds (uniformly) by the invariance
principles proved in \cite{SS, BB, MP}. 
 Assumption \ref{assu}(b) holds by Theorem 1.1 of \cite{B1}.

For $x \in \bZ^d$, let $S_x$ be the smallest integer $n$
such that 
$B_d(x,R)$ is very good with $N_{B_d(x,R)}^{2d+4}< R$
for all $R\ge n$.
(If $x \not\in \Ci$ we take $S_x=0$.)
Then by Theorem 2.18 and Lemma 2.19 of  \cite{B1} there
exists $\gamma=\gamma_d>0$
such that 
\begin{equation}\label{sxbnd}
\bP(S_x \ge n) \le c \exp(-c n^\gamma).
\end{equation}
In particular, we have that $S_x < \infty$ for all
$x \in \Ci$, $\bP$-a.s.
By Theorem \ref{hk2phi}, the PHI holds for
$Q(x,R,R^2)$ for all $R \ge S_x$, and Assumption \ref{assu}(c) holds. \\
(b) Assumption \ref{ass1}(e) holds 
by results in \cite{B1} -- see Proposition 2.17(d),
Lemma 2.19 and Remark 2 following Lemma 2.19.
\qed

In the results which follow, we have not made any effort
to obtain the best constant $\gamma$ in the various
bounds of the form $\exp(-n^\gamma)$.

\begin{lemma} \label{perc-hol}
With $\bP$-probability 1,
$\lim_{r \to \infty} h_\om(r) r^{-1/2} =0$, and so
Assumption \ref{assu}(d) holds.
\end{lemma}

\proof 
Let $M_0$ be the random variable given in Lemma 2.19
of \cite{B1}. Let $\alpha=1/4$, and note that 
$\beta=1-2(1+d)^{-1}>1/3$. Therefore
$$ \bP_0( M_0 \ge n) \le c \exp(-c n^{\alpha/3}), $$
and if $M_0\le n$ then the event
$D(Q,\alpha)$ defined in (2.21) of \cite{B1} holds
for every cube of side $n$ containing 0.
It follows from this
(see (2.20) and the definition of $R(Q)$ on p. 3040 in  \cite{B1})
that every cube of side greater than $n^\alpha$ in
$[-n/2,n/2]^d$ intersects $\Ci$. 
Thus
\begin{equation}\label{h-bound}
 \bP_0( h_\om(n) \ge n^\alpha) \le c  \exp(-c n^{\alpha/3}), 
\end{equation}
and using Borel-Cantelli we deduce that
$\lim_{r \to \infty} h_\om(r) r^{-1/2} =0$ $\bP_0$-a.s. \qed

\begin{lemma} \label{ass-f}
Let $x \in \bR^d$. 
With $\bP$-probability 1, Assumption \ref{ass1}(f)
holds.
\end{lemma}

\proof 
Let $F_n = \{  g^\om_n(x) \in \Lam_n(x,1) \}$,
and
$B_n = \{ S_{g^\om_n(x)} > n^{1/3} \}$.
If $F_n^c$ occurs, then a cube side $n$ containing $\Lam_n(x,1)$ 
has a hole greater than $n^{1/2}$. So, by \eqref{h-bound}
$$ \bP( F_m^c)  \le c e^{-c n^{1/3}}. $$
Let $Z_n = \max_{z \in \Lam_n(x,1)} S_x$.
Then
$$ B_n \subset  F_n^c \cup \{ Z_n >  n^{1/3} \}, $$
so using \eqref{sxbnd}
$$ \bP( B_n) \le
 c e^{-c' n^{1/3}} + c n^{d/2} e^{-c' n^{\gamma/3}}, $$
and by Borel-Cantelli Assumption 
\ref{ass1}(f) follows. \qed

\smallskip It remains to prove  Assumption \ref{ass1}(d). If instead
we wanted to control
$|\Lam_n|/(n^{1/2}\kappa)^d$
then we could use results in \cite{CM, DP}. Since the arguments
for $\mu(\Lam_n)$ are quite similar, we only give a sketch
of the proof.

\begin{lemma} \label{masseq}
Let $x \in \bR^d$. There exists $a>0$ such that with $\bP$-probability 1, 
 \begin{equation}\label{mulam}
 \frac{\mu(\Lam_n(x,r))}{ (2 n^{1/2}r)^d} \to a \qq
 \hbox{ as } n \to \infty,
\end{equation}
and so Assumption \ref{ass1}(d) holds.
\end{lemma}

\proof For a cube $Q \subset \bZ^d$ write $s(Q)$ for the length of the 
side of $Q$. Let $\pd_i Q= \pd(\bZ^d-Q)$ be the `internal boundary' of $Q$,
and $Q^0 = Q - \pd_i Q$. Recall that $\mu_x$ is the number of 
open bonds adjacent to $x$, and set
\begin{align*}
  M(Q) =\{ x \in Q^0: x \leftrightarrow \pd_i Q \}, \qq
  V(Q) = \mu( M(Q)).
\end{align*}
Note that if $x\in Q$ and $x$ is connected by an open path 
to $\pd_i Q$ then $x$ is connected to $\pd_i Q$ by an open path
inside $Q$. Thus the event $x \in M(Q)$ depends only on the
percolation process inside $Q$. So
if $Q_i$ are disjoint cubes, then the $V(Q_i)$ are
independent random variables. Let $\sC_k$ be a cube of side length $k$
and set
\begin{equation*}
 a_k = \bE k^{-d} V(C_k). 
\end{equation*}
By the ergodic theorem there exists $a$ such that,
$\bP$-a.s.,
\begin{equation}\label{elim}
 \lim_{R \to \infty} \frac{V(H(0,R/2))}{R^d}  \to a, \q 
\hbox{ $\bP$-a.s. and in $L^1$}.
\end{equation}
In particular, $a=\lim a_k$. 
Since $\Ci$ has positive density, it is clear that $a>0$.

We have
$$ \mu( Q \cap \Ci) \le V(Q) + c_1 s(Q)^{d-1}. $$
Let $\eps>0$. Choose $k$ large enough so that $c_1/k \le  \eps$,
and $a_k \le a + \eps$. 

Now let $Q$ be a cube of side $nk$, and let $Q_i$, $i=1, \dots n^d$
be a decomposition of $Q$ into disjoint sub-cubes each of side $k$.
Then
\begin{align*}
 (nk)^{-d} \mu(Q \cap \Ci) - a_k
&\le  (nk)^{-d} \sum_i \mu(Q_i \cap \Ci)  -a_k \\
&\le c_1 k^{-1} +  n^{-d} \sum_i (k^{-d} V(Q_i)- a_k).
\end{align*}
As this is a sum of i.i.d. mean 0 random variables, it follows that
there exists $c_2(k,\eps)>0$ such that
\begin{equation}\label{muub}
  \bP(  (nk)^{-d} \mu(Q \cap \Ci) > a +3 \eps) \le
 \exp ( -c_2(k,\eps) n^d).
\end{equation} 

The lower bound on $ \mu(Q \cap \Ci)$ requires a bit more
work. We call a cube $Q$ `$m$-good' if
the event $R(Q)$ given in \cite{AP} or p. 3040 of \cite{B1} holds, 
and 
$$ \mu(\Ci \cap Q) \ge (a-\eps) s(Q)^d. $$
Let $p_k$ be the probability a cube of side $k$ is $m$-good.
Then by (2.24) in \cite{AP}, and \eqref{elim},
$\lim p_k =1$.
As in \cite{AP} we can now divide $\bZ^d$ into disjoint macroscopic
cubes
$T_x$ of side $k$, and consider an associated site percolation 
process where a cube is occupied if it is $m$-good.
We write $\sC^*$ for the infinite cluster for this process.
Let $Q$ be a cube of side $nk$, and $T_x$ be the $n^d$ disjoint 
sub-cubes of side $k$ in $Q$. Then
\begin{equation}
 \mu( \Ci \cap Q) \ge \sum_x  \mu( \Ci \cap T_x)
 \ge (a-\eps) k^d \#\{x: T_x \in \sC^*, T_x \subset Q\}. 
\end{equation}
By Theorem 1.1 of \cite{DP} we can choose $k$ large enough
so there exists a constant $c_3(k,\eps)$ such that
\begin{equation}
\bP( n^{-d}  \#\{x: T_x \in \sC^*, T_x \subset Q\} < 1-\eps)
 \le \exp( - c_3(k,\eps) n^{d-1}).
\end{equation}
It follows that
\begin{equation}\label{mulb}
 \bP\big( (nk)^{-d}  \mu( \Ci \cap Q) < a-(1+a)\eps \big)
 \le  \exp( - c_3(k,\eps) n^{d-1}).
\end{equation}
Combining \eqref{muub} and  \eqref{mulb}, and using
Borel-Cantelli gives \eqref{mulam}. \qed

\ms\noindent  {\sl Proof of Theorem \ref{thm:llt-p}.}
By Lemmas \ref{easy-a}, \ref{ass-f} and \ref{masseq} 
Assumption \ref{ass1} holds for all $x \in \bQ^d$, $\bP$-a.s.,
and so also  $\bP_0$-a.s.
Therefore using Lemma \ref{easy-a} we have that
Assumption \ref{assu} holds   $\bP_0$-a.s., so \eqref{llt-pv}
follows from Theorem \ref{thm:unifllt}. \qed

\ms\noindent  {\sl Proof of Theorem \ref{thm:llt-pi}.}
The discrete time case  is given by Theorem  \ref{thm:llt-p}.
For continuous time, since Assumption \ref{assu} holds
 $\bP_0$-a.s.,
\eqref{llt-ctsver} follows from Theorem \ref{thm:unifllt-cts}. 
Since $a$ is given by \eqref{a-eqn}, and $\mu$ is the same
for $Y$ and the myopic walk, the constant $a$ in 
\eqref{llt-ctsver} is the same as for the myopic walk in 
\eqref{llt-pvi}. If $Z_t$ is a rate 1 Poisson process then we can
write $Y_t = X_{Z_t}$, and it is easy to check that the CLT for
$X$ implies one for $Y$ with the same diffusion constant $D$. \qed

\ms As a second application we consider the random 
conductance model in the case when the conductances are
bounded away from 0 and infinity.

Let $(\Omega, \sF, \bP)$ be a probability space. 
Let $K\ge 1$ and $\mu_e$, $e \in \bE_d$ 
be i.i.d.r.v. supported on $[K^{-1}, K]$.
Let also $\eta_x$, $x \in \bZ^d$ be i.i.d. random variables
on $[0,1]$, $F: R^{d+1} \to [K^{-1}, K]$, and
$\mu_{xx} = F(\eta_x, (\mu_{x \cdot}))$. For each
$\om \in \Omega$ let 
$X=(X_n ,n \in \bZ_+, P_\om^x, x \in \bZ^d)$ be the
SRW on $(\bZ^d, \mu)$ defined in Section 2, and
$p^\om_n(x,y)$ be its transition density.

\begin{theorem} \label{thm:llt-rcm}
Let $T_1>0$. Then there exist constants $a$, $D$ such that
$\bP_0$-a.s.,
\begin{equation}
\label{llt-rcm}
  \lim_{n\to\infty}  
\sup_{x\in \bR^d} \sup_{t \ge T_1} 
 |n^{d/2} \hat{p}^\om_{nt}(0,g^\om_n(x)) - 2 a^{-1} k^{(D)}_t(x)| =0.
\end{equation}
\end{theorem}

\proof As above, we just need to verify Assumption \ref{assu}.
The invariance principle in \cite{SS} implies the uniform CLT, giving
(a). Since $\mu_e$ are bounded away from 0 and infinity, the results
of \cite{Del} immediately give the PHI (with $S(x)=1$ for all $x$)
and heat kernel upper bound \eqref{fgub}, so giving 
Assumption \ref{assu}(b) and (c), as well as
Assumption \ref{ass1}(f).
As $\sG=\bZ^d$,  Assumption \ref{assu}(d) and
Assumption \ref{ass1}(e) hold.

It remains to verify Assumption \ref{ass1}(d), but this holds 
by an argument similar to that in Lemma \ref{masseq}. 
\qed

\section{Green's functions for percolation clusters}

We continue with the notation and hypotheses of Section
\ref{sec-perc}, but we take $d\ge 3$ throughout this section.
 The Green's function can be defined by
\begin{equation}
 g_\om(x,y) = \int_0^\infty q_t^\om(x,y) dt. 
\end{equation}
By Theorem \ref{disc}(c) 
$g_\om(x,y)$ is $\bP$-a.s. finite for all $x,y \in \Ci$.
We have that  $g_\om(x, \cdot)$ satisfies 
\begin{equation}\label{green-eqn}
  \sL g_\om(x, y) =
\begin{cases}
 0  &\hbox{ if } y \neq x , \\
 -1/\mu_x  & \hbox{ if } y=x. 
\end{cases}
\end{equation}
Since any bounded harmonic function is constant (see \cite{BLS}
or \cite[Theorem 4]{B1}), these equations have, $\bP$-a.s., a
unique solution such that $g_\om(x,y)\to 0$ as $|y|\to \infty$.
It is easy to check that the Green's function for the myopic and blind
ants satisfy the same equations, so the Green's function for the 
continuous time walk $Y$, and the myopic and blind ant discrete time
walks are the same.

We write $d_\om(x,y)$ for the graph distance on $\Ci$. 
By Lemma 1.1 and Theorem 1 of \cite{B1} there exist
$\eta>0$, constants $c_i$ and r.v. $T_x$ such that
\begin{equation}\label{txtail}
 \bP(T_x \ge n) \le c e^{-c_1 n^\eta},
\end{equation}
 so that the 
following bounds on $q^\om_t(x,y)$ hold:
\begin{align}\label{lrub}
 q_t(x,y) &\le  c_2 \exp(-c_3 d_\om(x,y) (1+ \log \fract{d_\om(x,y)}{t})), 
\q 1 \le t \le d_\om(x,y),\\
\label{cvub}
 q_t(x,y) &\le  c_4 e^{ -c_5 d_\om(x,y)^2/t}, \qq d_\om(x,y) \le t, \\
\label{gbq}
  c_6 t^{-d/2} e^{-c_7 |x-y|^2/t} &\le 
 q^\om_t(x,y) \le c_8 t^{-d/2} e^{-c_9 |x-y|^2/t}, \q t\ge  T_x \vee |y-x|.
\end{align}
We can and will assume that $T_x \ge 1$ for all $x$.

\begin{lemma} Let $x, y \in \Ci$, and $\delta \in (0,1)$. Then
\begin{align}\label{qtint-1}
  \int_0^{d_\om(x,y)} q_t^\om(x,y) dt &\le c_1 e^{-c_2 |x-y|},   \\
\label{qtint-2}
  \int_{d_\om(x,y)}^{T_x} q_t^\om(x,y) dt &\le c_3 T_x e^{-c_4 |x-y|^2/T_x}.
\end{align}
\end{lemma}

\proof 
Using \eqref{lrub} and \eqref{cvub} we have
\begin{align*} 
 \int_0^{d_\om(x,y)} q_t^\om(x,y) dt &\le
 \int_0^{d_\om(x,y)} c \exp(-c d_\om(x,y) ) dt \le c e^{-cd_\om(x,y)}, \\
  \int_{d_\om(x,y)}^{T_x} q_t^\om(x,y) dt &\le
   \int_{d_\om(x,y)}^{T_x} c e^{-c d_\om(x,y)^2 /t} dt 
\le c T_x e^{-c d_\om(x,y)^2/T_x},
\end{align*}
and since $d_\om(x,y) \ge c |x-y|$ this gives
\eqref{qtint-1} and  \eqref{qtint-2}.

\begin{proposition}\label{g-crude} Let $x,y \in \Ci$, with $x\neq y$. 
Then 
there exist constants $c_i$ such that
\begin{equation}\label{gcrudea}
   \frac{c_1}{|x-y|^{d-2}} \le
g_\om(x,y) \le \frac{c_2}{|x-y|^{d-2}} \qq 
\hbox{ if  } |x-y|^2 \ge  T_x (1 + c_3 \log |x-y|).
\end{equation}
Further, for $x, y \in \bZ^d$,
\begin{align}\label{gcrudeb}
   \frac{c_4}{1\vee |x-y|^{d-2}} \le
 &\bE \big( g_\om(x,y)| x, y \in \Ci \big) \le \frac{c_5}{1 \vee |x-y|^{d-2}},  \\
\label{gkbnd}
  &\bE \big( g_\om(x,x)^k| x \in \Ci \big) \le c_6(k).
\end{align}
\end{proposition}

\proof Note first that, by \eqref{gbq}
\begin{equation} \label{txub}
\int_{T_x}^\infty q_t^\om(x,y)dt 
\le \int_{0}^\infty c t^{-d/2} e^{-c |x-y|^2/t} dt \le c'|x-y|^{2-d}. 
\end{equation} 
Combining \eqref{qtint-1}, \eqref{qtint-2} and \eqref{txub} we obtain
\begin{equation}\label{g-ub-a}
 g_\om(x,y) \le c' e^{-c|x-y|} +   c T_x e^{-c_6|x-y|^2/T_x} + c|x-y|^{2-d}.
\end{equation}
Taking $c_3=d/c_6$ gives
$$  T_x e^{-c_6 |x-y|^2/T_x} \le c |x-y|^2 e^{- d \log |x-y|}
 \le c |x-y|^{2-d}, $$
and this gives the upper bound in \eqref{gcrudea}.
For the lower bound in \eqref{gcrudea} 
we note that since $T_x \le |x-y|^2$
\begin{align}\label{g-lower-a}
 g_\om(x,y) \ge \int_{|x-y|^2}^\infty  q_t^\om(x,y) dt
 \ge  \int_{|x-y|^2}^\infty   c t^{-d/2} e^{-c |x-y|^2/t} dt
 = c' |x-y|^{2-d} .
\end{align}

We now turn to \eqref{gcrudeb}.  Choose $k_0$ such that
$\bP( T_x \le k_0) \ge \half$. Then
\begin{align}\label{eg-lower-a}
 \bE^x g_\om(x,y) \ge 
 \bE^x \Big( \int_{T_x}^\infty q^\om_t(x,y)dt; T_x \le k_0 \Big)
 \ge \half \int_{k_0}^\infty c t^{-d/2} e^{-c |x-y|^2/t} dt. 
\end{align}
If $|x-y|^2 \ge k_0$, then the final term in \eqref{eg-lower-a}
is bounded below by $c |x-y|^{2-d}$ in the same way as in \eqref{eg-lower-a},
while when $|x-y|^2 \le k_0$ we have
\begin{equation}
  \bE^x g_\om(x,y) \ge c  \int_{k_0}^\infty c t^{-d/2} e^{-c |x-y|^2/t} dt
 \ge c e^{-c |x-y|^2/k_0} k_0^{1-d/2} \ge c',
\end{equation}
which gives the lower bound in \eqref{gcrudeb}.
For the averaged upper bound, note first that
\begin{equation}\label{gpower}
  g_\om(x,x) = \int_0^\infty q_t(x,x) dt \le c T_x + 
 \int_{T_x}^\infty c t^{-d/2} dt \le c' T_x.  
\end{equation} 
So for any $k\ge 1$, by \eqref{txtail}
$$ \bE( g_\om(x,x)^k| x \in \Ci) \le c(k) \bE( T_x^k | x \in \Ci)\le c'(k), $$
proving \eqref{gkbnd},
and (taking $k=1$) the upper bound in \eqref{gcrudeb} when $y=x$.
 
Now let $y\neq x$ and $F =\{ |x-y|^2 \le T_x(1+ c_{6.2.3} |x-y|) \}$. 
Then writing $\bE_{xy}(\cdot) = \bE(\cdot | x,y \in \Ci)$, and using
\eqref{gcrudea}, \eqref{gpower}, the fact that $g_\om(x,y) \le g_\om(x,x)$
and \eqref{txtail},
\begin{align*}
 \bE_{xy} g_\om(x,y) &= \bE_{xy}(  g_\om(x,y);F) +  \bE_{xy}(  g_\om(x,y);F^c) \\
  &\le c |x-y|^{2-d} + (\bE_{xy}(g_\om(x,y)^2))^{1/2} \bP_{xy}(F^c)^{1/2} \\
  &\le c |x-y|^{2-d} + (\bE_{xy}(g_\om(x,x)^2))^{1/2} c e^{-c |x-y|^{\eta/3}}
 \le c'  |x-y|^{2-d},
\end{align*}
proving \eqref{gcrudeb}.\qed

\ms To prove that $|y|^{d-2} g_\om(0,y)$
has a limit as $|y|\to \infty$ we use  Theorem \ref{thm:llt-pi}.
Write $ k_t(x) = k^{(D)}_t(x)$, where $D$ is the
constant in \eqref{llt-ctsver}.

\begin{lemma}\label{locallim}
Let $\eps>0$. Then for $\bP$-a.a. $\om \in \Omega_0$
there exists $a>0$ and $N=N(\eps,\om)$ such that
\begin{equation}\label{qk-bnd}
 | q^\om_t(0,y) - a^{-1} k_t(y) | \le \eps t^{-d/2} 
\q \hbox{ for all } t \ge N, \, y \in \Ci(\om). 
\end{equation}
\end{lemma}

\proof
By Theorem \ref{thm:llt-pi}.
there exists $N$ such that
\begin{equation} \label{llt-ctsver2}
\sup_{x\in \bR^d} \sup_{s \ge 1} 
 \Big|n^{d/2}q^\om_{ns}(0,g^\om_n(x))-  a^{-1} k_s(x)\Big| \le
 \eps \hbox{ for } n \ge N. 
\end{equation}
Let $n=N$, $s=t/n$ and $x = n^{-1/2} y$, so that
$g_n(x)=y$. Then noting that $k_s(x)= n^{d/2} k_t(y)$
\eqref{qk-bnd} follows. \qed

Let $|z|=1$ and 
\begin{equation}
  C = a^{-1} \int_0^\infty k_t(z) dt = (D a)^{-1} 
 \int_0^\infty (2\pi s)^{-d/2} e^{-1/2s}ds = \frac{\Gamma(\frac{d}{2}-1)}{2\pi^{d/2}aD}. 
\end{equation}

\sm {\it Proof of Theorem \ref{thm:green}.}
(a) This was proved as Proposition~\ref{g-crude}.

\sm (b) Let $\delta \in (0,1)$, to be chosen later.
For $y \in \Ci$ we set $t_1=t_1(y) = \delta |y|^2$, and
$t_2=t_2(y) = |y|^2/\delta $.
Then
\begin{equation}
 g_\om(0,y) = \int_0^{t_1}q^\om_t(0,y)dt + 
\int_{t_1}^{t_2} q^\om_t(0,y)dt + \int_{t_2}^\infty q^\om_t(0,y)dt
= I_1+I_2+I_3.
\end{equation}
As in Proposition \ref{g-crude} we have, 
using \eqref{qtint-1} and \eqref{qtint-2}, that 
provided $|y| \ge T_0$,
\begin{align}
 I_1 &\le c e^{-c|y|} + c T_0e^{-c |y|^2/T_0} + 
 \int_0^{\delta |y|^2} c t^{-d/2} e^{-c |y|^2/t}dt \\ 
 &\le   c e^{-c|y|} +c |y| e^{-c|y|} + 
c|y|^{2-d}  \int_0^{\delta } s^{-d/2} e^{-c_1/s}ds \\ 
 &\le c e^{-c|y|} + c|y|^{2-d} e^{-c_1/2\delta}.
\end{align}
Also
\begin{equation}
 I_3 \le \int_{|y|^2/\delta}^\infty c t^{-d/2} e^{-c |y|^2/t}dt
 = c \delta^{d/2-1} |y|^{2-d}.
\end{equation}
So there exist $M_1<\infty$ and $\delta>0$ so that
\begin{equation}
  I_1+I_3 \le \half \eps C |y|^{2-d} \qq \hbox{ when } |y|\ge M_1.
\end{equation}

Now let $\eps'>0$, and let $N=N(\eps')$ be given by Lemma \ref{locallim}.
For $I_2$ we have, provided $t_1\ge N$
\begin{align}\nonumber
 I_2 &\le \int_{t_1}^{t_2} (\eps' t^{-d/2} + a^{-1} k_t(y)) dt 
\le  c \eps' t_1^{1-d/2} + \int_{t_1}^{t_2} a^{-1} k_t(y) dt \\
&\le  c \eps' \delta^{1-d/2} |y|^{2-d} + C |y|^{2-d}. 
\end{align}
Taking 
$\eps'= \half (C/c) \eps \delta^{d/2-1}$ 
gives the upper bound in \eqref{greenlim2}. This bound
holds provided $|y| \ge M_1 \vee T_0 $ and
$\delta |y|^2 \ge N(\eps')$, Thus the upper bound in
 \eqref{greenlim2} holds provided
\begin{equation}\label{y-rest}
 |y| \ge T_0 \vee M_1 \vee (\delta^{-1} N(\eps'))^{1/2}. 
\end{equation}

For the lower bound, note that
\begin{equation}
 C|y|^{2-d}- \int_{t_1}^{t_2} a^{-1} k_t(y) dt 
\le c |y|^{2-d}( e^{-c/\delta} + \delta^{d/2-1} ). 
\end{equation}
So if \eqref{y-rest} holds then
\begin{align*}
 g_\om(0,y) \ge I_2 &\ge
 \int_{t_1}^{t_2} (- \eps' t^{-d/2} + k_t(y)) dt\\
 &\ge |y|^{2-d} \Big( C  -  c \eps' \delta^{1-d/2} 
 -  e^{-c/\delta} -\delta^{d/2-1} \Big),
 \end{align*}
proving the lower bound in \eqref{greenlim2}.

\sm (c) Let $\eps>0$, and $M$ be as in (a), and
$U_0= T_0(1+c_{6.2.3} \log |y|)$.
Then by Proposition \ref{g-crude}
\begin{align}\nonumber
 \bE_0 g_\om(0,y) &\le
 \bE_0(g_\om(0,y);M \le |y|) + \bE_0(g_\om(0,y); U_0 \le |y|< M) \\
 \nonumber
  & \qq \qq \qq \qq \qq \qq \qq\qq + \bE_0(g_\om(0,y); |y|< U_0) \\
\nonumber
 &\le \frac{(1+\eps)C}{|y|^{d-2}} + \frac{c_{6.2.2} }{|y|^{d-2}}\bP_0(M>|y|) 
+ (\bE_0 g_\om(0,y)^2)^{1/2} \bP_0( U_0 >|y|)^{1/2} \\
 \label{eg-tub}
 &\le \frac{(1+\eps)C + c_{6.2.2} \bP_0(M>|y|) }{|y|^{d-2}} 
 + c e^{- c |y|^{\eta/3}}
\end{align}
Also
\begin{equation}\label{eg-tlb}
 \bE_0 g_\om(0,y) \ge \bE_0(g_\om(0,y);M \le |y|)
 \ge \frac{(1-\eps)C}{|y|^{d-2}} \bP(M \le |y|).
\end{equation}
Combining \eqref{eg-tub} and \eqref{eg-tlb} completes the proof
of Theorem \ref{thm:green}. \qed

\appendix
\section{Appendix}\label{sec-app}

\sms In this appendix, we give a proof of the `balayage' formula
\eqref{bal-d1}-\eqref{bal-k1} used in the proof of the PHI in Section
\ref{sec-phi}.

Let $\Gam=(G,E)$ and $\mu$ be as in Section \ref{sec2}.
Let $B$ be a finite subset of $G$, and $B_1 \subset B$. 
Write $\ol B = B \cup \pd B$. Let $T \ge 1$, and
\begin{align*}
Q=(0,T]\times B, \q \ol Q = [0,T] \times \ol B, \q E=(0,T]\times B_1.
\end{align*}
\ms Set
\begin{equation}
  P^B f(x) = \sum_{y \in G} p^B_1(x,y) f(y)\mu_y, \q
   P f(x) = \sum_{y \in G} p_1(x,y) f(y)\mu_y,
\end{equation}
for any function $f$ on $G$

For a space-time function $w(r,y)$ we will sometimes write $w_r(y)=w(r,y)$.
Let
\begin{equation}
 Hw(n,x) = w(n,x) - P w_{n-1}(x).
\end{equation} 
Then $w$ is caloric in a space-time region $F\subset \bZ \times G$ 
if and only if $Hw(n,x)=0$ for $(n,x) \in F$.
Let $\sD$ be the set of non-negative functions $v(n,x)$ on $\ol Q$ such that
$v=0$ on $\ol Q-Q$ and $v$ is caloric on $Q-E$.
In particular we have $v(0,x)=0$ for $v \in \sD$.

\begin{lemma} \label{AL1}
Let $w(r,y)\ge 0$ on $\ol Q$, with $w=0$ on $\ol Q - E$, and let $v=v(n,x)$
be given by 
\begin{equation}\label{vA-def}
 v(n,x) = 
\begin{cases} 
 \sum_{r=1}^n P^B_{n-r} w_r(x), &\hbox{ if }  (n,x) \in Q \\
 0  &\hbox{ if }  (n,x) \not\in Q. 
\end{cases}
\end{equation}
Then $v \in \sD$, and 
\begin{equation}\label{Hv}
  H v(n,x) = w(n,x), \q (n,x) \in Q.
\end{equation}
\end{lemma}

\proof It is clear that $v\ge 0$, and that $v=0$ on $\ol Q-Q$.
If $x \in B$ then it easy to check that
$ P P^B_{m} f(x) =P^B_{m+1}f(x)$.
Let $(n,x)\in Q$, so $1\le n \le T$ and $x \in B$. Then
\begin{align}\nonumber
 Hv(n,x) 
 &= \sum_{r=1}^n P^B_{n-r} w_r(x) - P\big(\sum_{r=1}^{n-1} P^B_{n-1-r} w_r \big)(x)\\
 &=  \sum_{r=1}^n P^B_{n-r} w_r(x) - \sum_{r=1}^{n-1}  P^B_{n-r} w_r(x)
 =w_n(x).
\end{align}
This proves \eqref{Hv}, and as $w(n,x)=0$ when $x \in B-B_1$
we also deduce that $v$ is caloric in $Q-E$, proving that $v \in \sD$. \qed

\begin{lemma} \label{AL2}
 Let $u, v \in \sD$ satisfy 
$Hu(n,x)=Hv(n,x)$ for $(n,x) \in Q$. Then $u=v$ on $\ol Q$. 
\end{lemma}

\proof We have $u=v=0$ on $\ol Q -Q$.  We write
$u_k=u(k,\cdot)$. First note that $u_0=v_0$.  If $u_k=v_k$ and 
$x \in B$ then
$$ u(k+1,x)= Hu(k+1,x) + P u_k(x) = Hv(k+1,x) + P v_k(x), $$
so that $u_{k+1}=v_{k+1}$. \qed

Let $Z$ be the space-time process on $\bZ \times G$
given by $Z_n=(I_n, X_n)$, where $X$ is the SRW on $\Gamma$,
$I_n=I_0-n$, and $Z_0=(X_0,I_0)$ is the starting point of $Z$.
We write $\hat E^{(n,x)}$ for the law of $Z$ started
at $(n,x)$. 
Let $u(n,x)$ be non-negative and caloric on $Q$.
Then the r\'eduite $u_E$ is defined by 
\begin{equation} \label{red2}
 u_E(n,x)= \hat E^{(n,x)}\big( u(I_{T_E}, X_{T_E}); T_E< \tau_Q \big),   
\end{equation}
where
\begin{equation}
  T_E= \min\{k \ge 0: Z_k\in E \}, \q \tau_Q= \min\{k \ge 0: Z_k \not\in Q \}.
\end{equation}

\begin{lemma} \label{AL3} $u_E \in \sD$.
\end{lemma}

\proof If $(n,x) \in \ol Q -Q$ then $\hat P^{(n,x)}(\tau_Q=0)=1$,
so $u_E(n,x)=0$. It is clear from the definition \eqref{red2} 
that $u_E$ is caloric on $Q-E$, and that $u_E\ge 0$. \qed

\begin{proposition}\label{balay} Let $1\le n \le T$. Then 
\begin{equation}\label{bal-d1a}
 u_E(n,x) = \sum_{y \in B} \sum_{r=1}^n p^B_{n-r}(x,y) k(r,y) \mu_y,
\end{equation}
where
\begin{equation}\label{bal-k1a}
   k(r,y) = 
\begin{cases}
 \sum_{z \in B} p^B_1(y,z) (u(r-1,z)-u_E(r-1,z)) \mu_z, 
 &\hbox{ if } y \in B_1, \\
 0,  & \hbox{ if } y \in B-B_1. 
\end{cases}
\end{equation}
\end{proposition}

\proof
Let $k_r(y)=k(r,y)$ be defined by \eqref{bal-k1a} for $r \ge 1$. Set 
\begin{equation}
  v(n,x) =  \sum_{r=1}^n P^B_{n-r} k_r(x).
\end{equation} 
By Lemma \ref{AL1} we have $v \in \sD$. 
To prove that
$v=u_E$ it is sufficient, by Lemma \ref{AL2} to prove that
$Hv(n,x)=Hu_E(nx,)$ for $(n,x)\in Q$. 

We have $Hv(n,x) = k(n,x)$ on $Q$ by \eqref{Hv}. 
If $x \in B-B_1$ then $k(n,x)=0$, while since $u_E$ is caloric
in $Q-E$ we have $Hu_E(n,x)=0$. If $x \in B_1$ then
as $u=u_E$ on $E$, and $u$ is caloric on $Q$,
\begin{align*}
  Hu_E(n,x)&= u_E(n,x) - P u_E(n-1,x) \\
        &=  u(n,x) - P u_E(n-1,x) = P u(n-1,x) - P u_E(n-1,x) \\
  &= P^B_1 (u-u_E)(n-1,x).
\end{align*}
So we deduce that $v=u_E$. \qed

\sms If $y \in B_1$ then
the $r=1$ term of \eqref{bal-d1a} can be written
\begin{equation}
 \sum_{y\in B} p^B_{n-1}(x,y)\mu_y  (\sum_{z\in B} p^B_1(y,z)\mu_x u(0,z))
 =  \sum_{z\in B} \mu_x u(0,z) p^B_{n}(x,z), 
\end{equation}
so that \eqref{bal-d1a} can be rewritten as
\begin{equation}\label{bal-d1b}
 u_E(n,x) =  \sum_{y \in B} p^B_n(x,y) u(0,y)\mu_y 
  +\sum_{y \in B} \sum_{r=2}^n p^B_{n-r}(x,y) k(r,y) \mu_y,
\end{equation}
which is the form given in \eqref{bal-d1}.

\med
{\bf Acknowledgement} We are grateful to R. Cerf for remarks 
concerning Lemma \ref{masseq} and to J. Cerny for asking about
Green's functions.

\vskip 0.5 truein

\noindent MB: Department of Mathematics,
University of British Columbia,
Vancouver, BC V6T 1Z2, Canada.

\noindent BH: Mathematical Institute, University of Oxford, 24-29 St
Giles, Oxford OX1 3LB, UK.

\end{document}